\documentclass{article}
\usepackage{amssymb,amsmath,bm,geometry,graphics,graphicx,color,url}

\def\no{\noindent}
\def\pmatrix{\left(\begin{array}}
\def\endpmatrix{\end{array}\right)}
\def\RR{\mathbb{R}}

\def\I{{\cal I}}

\def\M{{\cal M}}

\def\P{{\cal P}}

\def\dd{\mathrm{d}}

\def\diag{\mathrm{diag}}

\newtheorem{theo}{Theorem}
\newtheorem{lem}{Lemma}
\newtheorem{cor}{Corollary}
\newtheorem{rem}{Remark}

\def\proof{\noindent\underline{Proof}\quad}
\def\QED{\mbox{~$\Box{~}$}}

\def\bq{u}
\def\bp{v}
\def\blam{{\lambda_0}}

\def\tq{{\tilde q}}
\def\tp{{\tilde p}}
\def\tlam{{\tilde\lambda}}

\def\bfgamma{{\bm{\gamma}}}

\def\bfpsi{{\bm{\psi}}}
\def\bfrho{{\bm{\rho}}}

\def\hb{\hat{b}}
\def\hc{\hat{c}}

\def\hH{\hat{H}}

\title{Line Integral Solution of Hamiltonian Systems with Holonomic Constraints}

\author{Luigi Brugnano\,$^a$,~ Gianmarco Gurioli\,$^a$,~Felice Iavernaro\,$^b$,~ Ewa B.\,Weinm\"uller\,$^c$ \\
~\\
$^a$ {\small Dipartimento di Matematica e Informatica ``U.\,Dini'', Universit\`a di Firenze}\\
{\small  Viale Morgagni 67/a, I-50134 Firenze, Italy.}\\
$^b$ {\small Dipartimento di Matematica, Universit\`a di Bari}\\
{\small Via Orabona 4, I-70125 Bari, Italy.}\\
$^c$ {\small Institute for Analysis and Scientific Computing,
Vienna University of Technology}\\
{\small A-1040 Wien, Austria.} }

\date{\normalsize\em -- Dedicated to John Butcher, on the occasion of his 84-th birthday --}

\begin{document}

\maketitle

\begin{abstract} In this paper, we propose a second-order energy-conserving approximation procedure for Hamiltonian systems with holonomic constraints.  The derivation of the procedure relies on the use of the so-called {\em line integral framework}. We provide numerical experiments to illustrate  theoretical findings.

\medskip
\no{\bf Keywords:} constrained Hamiltonian systems; holonomic constraints; energy-conserving methods; line integral methods; Hamiltonian Boundary Value Methods; HBVMs.\\

\smallskip
\no{\bf MSC:} 65P10, 65L80, 65L06.
\end{abstract}

\section{Introduction}\label{intro}
We consider the numerical approximation of a constrained Hamiltonian dynamics, described by the separable Hamiltonian
\begin{equation}\label{Hqp}
H(q,p) =  \frac{1}2 p^\top  M^{-1} p +U(q), \qquad q,p\in\RR^m,
\end{equation}
where $M$ is a symmetric and positive-definite matrix. The problem is completed by $\nu$ {\em holonomic} constraints,
\begin{equation}\label{gq}
g(q) = 0\in\RR^\nu,
\end{equation}
where we assume that $\nu< m$ holds. Moreover, we also assume that all points are regular for the constraints, i.e., $\nabla g(q)\in\RR^{m\times \nu}$ has full column rank or, equivalently, $\nabla g(q)^\top M^{-1}\nabla g(q)$ is nonsingular. For simplicity, both $U$ and $g$ are assumed to be analytic.

 It is well-known that the problem defined by (\ref{Hqp})--(\ref{gq}) can be cast in Hamiltonian form by defining the augmented Hamiltonian
\begin{equation}\label{Hqpl}
\hH(q,p,\lambda) = H(q,p) + \lambda^\top g(q),
\end{equation}
where $\lambda$ is the vector of Lagrange multipliers. The resulting constrained Hamiltonian system reads: \begin{equation}\label{constrp}
\dot q =  M^{-1}p, \qquad \dot p = -\nabla U(q)-\nabla g(q)\lambda, \qquad g(q)=0, \qquad  t\in[0,T],
\end{equation}
and is subject to consistent initial conditions,
\begin{equation}\label{q0p0}
q(0)=q_0,\qquad p(0)=p_0,
\end{equation}
such that
\begin{equation}\label{consist0}
g(q_0)=0, \qquad \nabla g(q_0)^\top  M^{-1} p_0 = 0.
\end{equation}
Note that the condition $g(q_0)=0$ ensures that $q_0$ belongs to the manifold
\begin{equation}\label{M}
\M = \left\{ q\in\RR^m: g(q)=0\right\},
\end{equation}
as required by the constraints, whereas the condition $\nabla g(q_0)^\top  M^{-1}p_0$ means that the motion initially stays on the tangent space to $\M$ at $q_0$. On a continuous level, this condition is satisfied by all points on the solution trajectory, since, in order for the constraints to be conserved,  
\begin{equation}\label{Dgp0}
\dot g(q) = \nabla g(q)^\top \dot q = \nabla g(q)^\top  M^{-1} p = 0,
\end{equation}
holds. These latter constraints are sometimes referred to as {\em hidden constraints}.

We stress that the condition (\ref{Dgp0}) can be conveniently relaxed for the numerical approximation. There, we only ask for $\nabla g(q)^\top  M^{-1}p$ to be suitably small along the numerical solution. Consequently, when solving the problem on the interval $[0,h]$, we require that the approximations,
\begin{equation}\label{q1p1}
q_1\approx q(h), \qquad p_1\approx p(h),
\end{equation} 
satisfy the conservation of both, the Hamiltonian and the constraints,
\begin{equation}\label{Hgcons}
H(q_1,p_1) = H(q_0,p_0), \qquad g(q_1) = g(q_0) = 0,
\end{equation}
and that the hidden constraints are relaxed to
\begin{equation}\label{nablagp}
\nabla g(q_1)^\top  M^{-1}p_1 = O(h^2).
\end{equation}
We recall that a formal expression for the vector $\lambda$ is obtained by an additional differentiating of (\ref{Dgp0}), i.e.,
\begin{equation}\label{lam0}
\ddot g(q) = \nabla^2 g(q) (M^{-1}p,M^{-1}p) - \nabla g(q)^\top M^{-1} \left[ \nabla U(q) + \nabla g(q)\lambda\right].
\end{equation}
Imposing the vanishing of this derivative yields 
\begin{equation}\label{lambda1}
\left[\nabla g(q)^\top M^{-1} \nabla g(q)\right]\lambda =  \nabla^2 g(q) (M^{-1}p,M^{-1}p) ~-~ \nabla g(q)^\top M^{-1} \nabla U(q).
\end{equation}
Consequently, the following result follows. 

\begin{theo}\label{lamex} The vector $\lambda$ exists and is uniquely determined, provided that the matrix $$\nabla g(q)^\top M^{-1} \nabla g(q)$$ is nonsingular.
In fact, in such a case, from (\ref{lambda1}), we obtain
\begin{eqnarray}\nonumber
\lambda &=&  \left[\nabla g(q)^\top M^{-1} \nabla g(q)\right]^{-1}\left[ \nabla^2 g(q) (M^{-1}p,M^{-1}p) ~-~ \nabla g(q)^\top M^{-1} \nabla U(q)\right]
\\[1mm]
&=&: \lambda(q,p). \label{lqp}
\end{eqnarray}
Note that, for later use, we have introduced the notation $\lambda(q(t),p(t))$ in place of $\lambda(t)$, to explicitly underline the dependence of the Lagrange multiplier on the state variables $q$ and $p$ at time $t$. 
\end{theo}

\begin{rem} We observe that an additional differentiation of (\ref{lambda1}) provides a differential equation for the Lagrange multipliers, which can be solved together with the original problem. However, this procedure is cumbersome in general, since it requires the evaluation of higher order tensors.
\end{rem}

\medskip
Numerical solution of Hamiltonian problems with holonomic constraints has been for a long time in the focus of interest. Many different approaches have been proposed such as the basic Shake-Rattle method \cite{Shake1977,Rattle1983}, which has been shown to be symplectic \cite{LeSk1994}, higher order methods obtained via symplectic PRK methods \cite{Jay1996}, composition methods \cite{Reich1996,Reich1997}, symmetric LMFs \cite{CHL2013}. Further methods are based on discrete derivatives \cite{Gonzalez1999}, local parametrizations of the manifold containing the solution \cite{BCF2001}, or on projection techniques \cite{Seiler1998,WDLW2016}. See also \cite{LeRe1994,Seiler1999,Hairer2003,LBS2004} and the monographs \cite{BCP1996,LeRe2004,HLW2006,HaWa2010}.

In this paper we pursue a different approach, utilizing the so-called {\em line integral}, which has already been used when deriving the energy-conserving Runge-Kutta methods, for unconstrained Hamiltonian systems, cf.\ Hamiltonian Boundary Value Methods (HBVMs) \cite{BIT2009,BIT2010,BIT2011,BIT2012_2,BIT2015} and the recent monograph \cite{LIMbook2016}. Such methods have also been applied in a number of applications \cite{BCMr2012,BrIa2012,BIT2012,BIT2012_1,ABI2015,BFCI2015,BBGCI17}, and here are used to cope with the constrained problem (\ref{Hqp})--(\ref{gq}). Roughly speaking, the conservation of the invariant will be guaranteed by requiring that a suitable line integral vanishes. This line integral represents a discrete-time version of (\ref{Dgp0}). In fact, if we fix a stepsize $h>0$, then the conservation of the constraints (\ref{gq}) at $h$, starting from the point $q_0$ defined in (\ref{q0p0}), can be recast into 
$$g(q(h)) - \underbrace{ g(q(0)) }_{=0} = \int_0^h \nabla g(q(t))^\top\dot q(t)\dd t    ~=~ 0.$$
For the continuous solution, this integral vanishes since the integrand is identically zero due to (\ref{constrp}) and (\ref{Dgp0}).  However, we can relax this requirement in the context of a numerical method describing a {\it discrete-time} dynamics. In such a case, the conservation properties have to be satisfied only on a set of discrete times which are multiples of the stepsize $h$. Consequently, we consider a local approximation to $q(t)$, say $\bq(t)$, such that
$$\bq(0) = q_0, \qquad \bq(h) =: q_1\approx q(h),$$
and
$$g(q_1)-g(q_0) \equiv g(\bq(h))-g(\bq(0)) = \int_0^h \nabla g(\bq(t))^\top \dot\bq(t)\dd t = 0,$$
without requiring the integrand to be identically zero. This, in turn, enables a proper choice of the vector of the multipliers $\lambda$. As a result, we eventually obtain suitably modified HBVMs which enable to conserve both, the Hamiltonian and the constraints. We stress that the available efficient implementation of the original methods (see, e.g., \cite{BIT2011,BFCI2014,LIMbook2016}), which proved to be reliable and robust in the numerical solution of the unconstrained Hamiltonian problems, can now be adapted for dealing with the holonomic constraints.\\
%&&&&&&&&&&&&&&&&&&&&&&&&&&&&&&&&&&&&
The paper is organized as follows. In Section~\ref{polapp}, we provide the framework for devising the method via a suitable choice of the vector $\lambda$ of the Lagrange multipliers, which we approximate by a piecewise-constant function. In Section~\ref{poly} further simplification towards numerical procedure is discussed. Then, in Section~\ref{HBVMs}, we present a fully discrete method, resulting in a suitable modification of the original HBVMs. In Section~\ref{numtest}, numerical experiments are shown to illustrate how the method works for a number of constrained Hamiltonian problems. Section~\ref{end} contains the conclusions and possible future investigations.

\section{Piecewise-constant approximation of $\lambda$}\label{polapp}

In this section, we show that we can approximate the solution of problem (\ref{constrp}) on the interval $[0,h]$, $h=T/N$, by looking for a constant vector $\lambda\in\RR^{\nu}$ such that (\ref{Hgcons}) is satisfied. This is equivalent to require 
\begin{equation}\label{Hlgcons}
\hH(q_1,p_1,\lambda) = \hH(q_0,p_0,\lambda), \qquad g(q_1) = g(q_0) = 0,
\end{equation}
where the constant parameter $\lambda$ is chosen in such a way that the constraints $g(q_1)=0$ hold. We will show that this procedure provides us with a second order approximation of the original problem, which becomes exact when the true multiplier is constant. Consequently, we approximate the problem (\ref{constrp})--(\ref{q0p0}), by the local problem
\begin{equation}\label{constrp_lc}
\dot \bq =  M^{-1}\bp, \qquad \dot \bp = -\nabla U(\bq)-\nabla g(\bq)\blam, \qquad t\in[0,h],
\end{equation}
subject to the initial conditions, cf.\ (\ref{q0p0}),
\begin{equation}\label{u0v0}
\bq(0) = q_0, \qquad \bp(0) = p_0,
\end{equation}
satisfying (\ref{consist0}). By setting
\begin{equation}\label{q1p1def}
q_1:=u(h), \qquad p_1 :=v(h),
\end{equation}
the constant parameter $\blam$ is chosen to guarantee the conservation of the Hamiltonian and the constraints, i.e., (\ref{Hgcons}).
Starting from (\ref{q1p1def}), the procedure is then repeated on $[h,2h]$ and the following intervals. The convergence result is now formulated in the following theorem.  

\begin{theo}\label{capdcazz}
For all sufficiently small stepsizes $h>0$, the above procedure defines a sequence of approximations $(q_n,p_n)$ such that,
for all $n=1,2,\ldots:$
\begin{equation}\label{ordqp}
q_n = q(nh) + O(h^2), \qquad p_n = p(nh) + O(h^2), \qquad g(q_n)=0, \qquad \nabla g(q_n)^\top M^{-1}p_n = O(h^2).
\end{equation}
Moreover, $(q_{n+1},p_{n+1})$ is obtained from $(q_n,p_n)$ using a constant vector $\lambda_n$ such that
\begin{equation}\label{ordl}
\lambda_n = \lambda( q(nh), p(nh) ) + O(h),
\end{equation}
where $\lambda(q,p)$ is defined in (\ref{lqp}) and, consequently, ~$H(q_{n+1},p_{n+1}) = H(q_n,p_n).$
\end{theo}

\medskip
The aim of this section is to show (\ref{ordqp})--(\ref{ordl}). Let us first consider 
 the orthonormal basis on $[0,1]$ given by the shifted and scaled Legendre polynomials $\{P_j\}$,
\begin{equation}\label{orto}
P_j\in\Pi_j, \qquad \int_0^1 P_i(c)P_j(c)\dd c=\delta_{ij}, \qquad \forall i,j=0,1,\dots,
\end{equation}
along with the expansions,
\begin{eqnarray}\nonumber
 M^{-1}\bp(ch) &=& \sum_{j\ge0} P_j(c)\gamma_j(\bp), \qquad \nabla U(\bq(ch)) ~=~ \sum_{j\ge0} P_j(c)\psi_j(\bq),\\ \label{expHg}
\nabla g(\bq(ch)) &=& \sum_{j\ge0} P_j(c)\rho_j(\bq), \qquad c\in[0,1],
\end{eqnarray}
with
\begin{eqnarray}\nonumber
\gamma_j(\bp) &=&  M^{-1}\int_0^1 P_j(c)\bp(ch)\dd c, \qquad \psi_j(\bq) ~=~ \int_0^1 P_j(c)\nabla U(\bq(ch))\dd c,\\  \label{cgj}
\rho_j(\bq) &=& \int_0^1 P_j(c)\nabla g(\bq(ch))\dd c, \qquad j\ge 0.
\end{eqnarray}
Consequently, following the approach defined in \cite{BIT2012_2}, the differential equations in (\ref{constrp_lc}) can be rewritten as
\begin{equation}\label{constrp1}
\dot \bq(ch) = \sum_{j\ge0} P_j(c) \gamma_j(\bp), \qquad \dot \bp(ch) = -\sum_{j\ge0} P_j(c)[ \psi_j(\bq) +\rho_j(\bq)\blam], \qquad  c\in[0,1].
\end{equation}
Moreover, using the initial conditions (\ref{q0p0}), we formally obtain
\begin{equation}\label{qpch}
\bq(ch) = q_0 + h\sum_{j\ge0} \int_0^cP_j(x)\dd x\, \gamma_j(\bp), \quad  \bp(ch) = p_0 -h\sum_{j\ge0} \int_0^cP_j(x)\dd x[ \psi_j(\bq) +\rho_j(\bq)\blam], \quad  c\in[0,1].
\end{equation}
The following result is now cited from \cite[Lemma\,1]{BIT2012_2}.
\begin{lem}\label{Ohj} Let $G:[0,h] \rightarrow V$, with $V$ a vector space, admit a Taylor expansion at $0$. Then
$$\int_0^1 P_j(c) G(ch)\dd c = O(h^j), \qquad j\ge0.$$\end{lem}

As a straightforward consequence, one has the following result.

\begin{cor}\label{cor1} All coefficients specified in (\ref{cgj}) are $O(h^j)$. \end{cor}

Concerning the conservation properties of the approximations, the following result holds.

\begin{theo}\label{hH_cons} For all $\blam\in\RR^\nu$, the solution of (\ref{constrp_lc})--(\ref{q1p1def}) satisfies $$\hH(q_1,p_1,\blam) = \hH(q_0,p_0,\blam).$$\end{theo}
\proof For any given $\blam\in\RR^\nu$, it follows from (\ref{Hqpl}), (\ref{constrp_lc}), and (\ref{q1p1def}),
\begin{eqnarray*}
\lefteqn{
\hH(q_1,p_1,\blam) - \hH(q_0,p_0,\blam) ~=~ \hH(\bq(h),\bp(h),\blam) - \hH(\bq(0),\bp(0),\blam)}\\[2mm]
&=& \int_0^h \frac{\dd}{\dd t} \hH(\bq(t),\bp(t),\blam)\dd t~=~\int_0^h \left\{\hH_q(\bq(t),\bp(t),\blam)^\top\dot\bq(t) + \hH_p(\bq(t),\bp(t),\blam)^\top\dot\bp(t)\right\}\dd t \\
&=& h\int_0^1 \left\{[\nabla U(\bq(ch)) + \nabla g(\bq(ch))\blam]^\top\dot\bq(ch) + [ M^{-1}\bp(ch)]^\top\dot\bp(ch)\right\}\dd c\\
&=& h\int_0^1 \left\{[\nabla U(\bq(ch)) + \nabla g(\bq(ch))\blam]^\top\sum_{j\ge0}P_j(c)\gamma_j(\bp) ~-~\right.\\
&&\qquad \qquad \left. [ M^{-1}\bp(ch)]^\top\sum_{j\ge0}P_j(c) [ \psi_j(\bq)+\rho_j(\bq)\blam]\right\}\dd c\\
\end{eqnarray*}\begin{eqnarray*} %%%%%%%%%%%%%%%%%%%%%%%%%%%%%%%%%%%%%%%%%%%%%%%%%%% to split the formula
&=&h\sum_{j\ge0}\left\{\left( \int_0^1 P_j(c)[\nabla U(\bq(ch)) + \nabla g(\bq(ch))\blam]\dd c \right)^\top\gamma_j(\bp) ~-~\right.\\
&&\qquad \left. \left( M^{-1}\int_0^1 P_j(c)\bp(ch)\dd c\right)^\top[ \psi_j(\bq)+\rho_j(\bq)\blam]\right\}\\
&=&h\sum_{j\ge0}\left\{ [\psi_j(\bq)+\rho_j(\bq)\blam]^\top \gamma_j(\bp) -
\gamma_j(\bp)^\top [\psi_j(\bq)+\rho_j(\bq)\blam]\right\} = 0.{\hfil \QED}
\end{eqnarray*}
As observed above, the conservation of the Hamiltonian (\ref{Hqp}) is guaranteed, once the constraints are satisfied, i.e., $g(q_1)=0$. We now apply a {\em line integral} technique to determine the vector $\blam$ and formulate the following result describing the very first step of the approximation procedure.

\begin{theo}\label{lambda} Let us consider the problem (\ref{constrp_lc})--(\ref{u0v0}) and assume that $(q_0,p_0)$ is given such that,
\begin{itemize}
\item $\nabla g(q_0)^\top  M^{-1} \nabla g(q_0) \in\RR^{\nu\times \nu}$ is nonsingular;
\item $g(q_0) = 0$;
\item $\nabla g(q_0)^\top M^{-1} p_0 = 0$.
\end{itemize}
Then, for all sufficiently small $h>0$,  $\exists!\blam\in\RR^\nu$ such that the approximations in (\ref{q1p1def}) satisfy
\begin{itemize}
\item $g(q_1)=0$ and, therefore, $H(q_1,p_1)=H(q_0,p_0)$;
\item $\lambda_0 = \lambda(q_0,p_0) + O(h)$;
\item $q_1-q(h) = O(h^2)$, $p_1-p(h)=O(h^2)$;
\item $\nabla g(q_1)^\top M^{-1} p_1 = O(h^2)$.
\end{itemize}
\end{theo}

\begin{rem} Clearly, Theorem~\ref{lambda} is the discrete counterpart of Theorem~\ref{lamex}.
\end{rem}

Before showing Theorem~\ref{lambda}, we have to state the following preliminary results.

\begin{lem}\label{intP} Let us consider the polynomial basis (\ref{orto}). Then, we have 
\begin{equation}\label{intP1}
\int_0^1P_j(c)\int_0^c P_i(x)\dd x\,\dd c = \left( X_s \right)_{j+1,i+1}, \qquad i,j=0,\dots,s-1,
\end{equation}
where $\left( X_s \right)_{j+1,i+1}$ is the $(j+1,i+1)$ entry of the matrix
\begin{equation}\label{Xs}
X_s:=\pmatrix{rrrr}
\xi_0 & -\xi_1\\
\xi_1 &     0   &\ddots\\
         &   \ddots &\ddots & -\xi_{s-1}\\
         &              &\xi_{s-1} &0\endpmatrix, \qquad \xi_j = \frac{1}{2\sqrt{|4j^2-1|}},\quad j=0,\dots,s-1.
         \end{equation}
\end{lem}
\proof Since the integrand on the left-hand side in (\ref{intP1}) is a polynomial of degree at most $2s-1$, the integral can be computed exactly via the Gauss-Legendre formula of order $2s$. Let $c_1,\dots,c_s$ be the zeros of $P_s$ and $b_1,\dots,b_s$ be the corresponding weights. Then, introducing the matrices
\begin{equation}\label{PsIsOm}
\P_s = \left( P_{j-1}(c_i) \right), ~ \I_s = \left( \int_0^{c_i} P_{j-1}(x)\dd x\right), ~\Omega = \diag(b_1,\dots,b_s) ~\in\RR^{s\times s},
\end{equation}
and setting $e_i\in\RR^s$, the $i$-th unit vector, we have
$$
\int_0^1P_j(c)\int_0^c P_i(x)\dd x\,\dd c ~=~ \sum_{\ell=1}^s b_\ell P_j(c_\ell)\int_0^{c_\ell} P_i(x)\dd x ~\equiv~ e_{j+1}^\top \P_s^\top \Omega \I_s e_{i+1}.
$$
The result follows by observing that, due to the properties of Legendre polynomials \cite[Section\,1.4.3]{LIMbook2016}, $$\I_s = \P_s X_s, \quad \P_s^\top\Omega\P_s=I_s$$ follows, where $X_s$ is the matrix defined in (\ref{Xs}), and $I_s\in\RR^{s\times s}$ is the identity matrix. This yields
$$ e_{j+1}^\top \P_s^\top \Omega \I_s e_{i+1} = e_{j+1}^\top \P_s^\top \Omega \P_s X_s e_{i+1} = e_{j+1}^\top X_s e_{i+1}.\,\QED$$

We also need the following expansions.

\medskip
\begin{lem}\label{fanculo} From (\ref{constrp_lc}) and (\ref{cgj}), we conclude 
\begin{eqnarray*}
\rho_0(\bq) &=& \nabla g(q_0) + \frac{h}2\nabla^2 g(q_0) M^{-1}p_0 + O(h^2),\\[2mm]
\rho_0(\bq)^\top M^{-1}p_0 &=& \nabla g(q_0)^\top M^{-1} p_0 + \frac{h}2\nabla^2 g(q_0) (M^{-1} p_0,M^{-1}p_0) + O(h^2).
\end{eqnarray*}
\end{lem}
\proof The first equality follows from Lemma~\ref{Ohj} and Corollary~\ref{cor1},
\begin{eqnarray*}
\rho_0(\bq) &=& \int_0^1 \nabla g (\bq(ch))\dd c ~=~ \int_0^1\left[\nabla g(\bq(0)) + ch \nabla^2 g(\bq(0))\dot\bq(0) + O((ch)^2)\right]\dd c \\
&=& \nabla g(\bq(0)) + \frac{h}2 \nabla^2 g(\bq(0))\dot\bq(0) +O(h^2)\\
&=& \nabla g(q_0) + \frac{h}2\nabla^2 g(q_0)\sum_{j\ge0}  P_j(0)\gamma_j(\bp)  +O(h^2)\\
%\end{eqnarray*}\begin{eqnarray*} %%%%%%%%%%%%%%%%%%%%%%%%%%%%%%%%%%%%%%%%%%%%%%%%%%% to split the formula
&=& \nabla g(q_0) + \frac{h}2\nabla^2 g(q_0)\gamma_0(\bp) +O(h^2)\\[2mm]
&=& \nabla g(q_0) + \frac{h}2\nabla^2 g(q_0)M^{-1} \left[ p_0 + O(h)\right] +O(h^2)\\[2mm]
&=& \nabla g(q_0) + \frac{h}2\nabla^2 g(q_0)M^{-1}p_0 +O(h^2).
\end{eqnarray*}
The second statement follows by transposition and multiplication from the right by $M^{-1}p_0$.\,\QED
\bigskip

We now show the results formulated in Theorem~\ref{lambda}.
\bigskip

\proof (of Theorem~\ref{lambda}). ~ Let us assume that  $g(q_0)=0$ holds. Then, it follows from (\ref{constrp_lc})--(\ref{q1p1def}), 
\begin{eqnarray*}
\lefteqn{g(q_1) ~=~ g(q_1)- g(q_0) ~=~g(\bq(h)) - g(\bq(0)) ~=~ \int_0^h \frac{\dd}{\dd t} g(\bq(t))\dd t}\\
&=& \int_0^h \nabla g(\bq(t))^\top\dot\bq(t)\dd t ~=~ h\int_0^1 \nabla g(\bq(ch))^\top\dot\bq(ch)\dd c ~=~
h\int_0^1 \nabla g(\bq(ch))^\top\sum_{j\ge0} P_j(c) \gamma_j(\bp)\dd c\\
&=& h\sum_{j\ge0} \rho_j(\bq)^\top\gamma_j(\bp) ~=~
h\sum_{j\ge0} \rho_j(\bq)^\top  M^{-1}\int_0^1 P_j(c)\bp(ch)\dd c\\
&=& h\sum_{j\ge0} \rho_j(\bq)^\top  M^{-1} \int_0^1 P_j(c) \left\{p_0-h\sum_{i\ge0} \int_0^c P_i(x)\dd x[ \psi_i(\bq) +\rho_i(\bq)\blam]\right\}\dd c\\
%\end{eqnarray*}\begin{eqnarray*} %%%%%%%%%%%%%%%%%%%%%%%%%%%%%%%%%%%%%%%%%%%%%%%%%%% to split the formula
&=& h\sum_{j\ge0} \rho_j(\bq)^\top  M^{-1}p_0 \int_0^1 P_j(c) \dd c\\
&&~ -~h^2\sum_{i,j\ge0} \rho_j(\bq)^\top  M^{-1}[ \psi_i(\bq) +\rho_i(\bq)\blam] \int_0^1 P_j(c)\int_0^c P_i(x)\dd x\,\dd c.
\end{eqnarray*}
Due to (\ref{orto}),
$$\int_0^1 P_j(c)\dd c = \delta_{j0},$$ 
and according to (\ref{intP1})--(\ref{Xs}), we conclude

\begin{eqnarray}\nonumber
\lefteqn{g(\bq(h)) - g(\bq(0))~=~h\rho_0(\bq)^\top  M^{-1} \left\{ p_0 -h[ \xi_0(\psi_0(\bq)+\rho_0(\bq)\blam)
- \xi_1(\psi_1(\bq)+\rho_1(\bq)\blam)] \right\} }\\ \nonumber
&&~-~h^2\sum_{j\ge1}\rho_j(\bq)^\top  M^{-1}\left\{[
\xi_j(\psi_{j-1}(\bq)+\rho_{j-1}(\bq)\blam) - \xi_{j+1}(\psi_{j+1}(\bq)+\rho_{j+1}(\bq)\blam)]\right\} \\ \label{hGamma}
%&&~-~h^2\rho_{s-1}(\bq)^\top  M^{-1}[\xi_{s-1}(\psi_{s-2}(\bq)+\rho_{s-2}(\bq)\blam) ] ~=:~
&=:&\hat\Gamma(\bq,\bp,\blam,h) .
\end{eqnarray}
By virtue of (\ref{cgj}) and Corollary~\ref{Ohj}, 
\begin{eqnarray}\nonumber
\lefteqn{\frac{g(\bq(h))-g(\bq(0)) - h\rho_0(\bq)^\top M^{-1}p_0}{h^2} =}\\  \label{consist}
&& -\frac{1}2\left\{ \left[\rho_0(\bq)^\top M^{-1}\rho_0(\bq) + O(h)\right]\blam ~+~ \rho_0(\bq)^\top M^{-1}\psi_0(\bq) +O(h)\right\}
\end{eqnarray}
follows. Now, from (\ref{cgj}) and Lemma~\ref{fanculo}, we have 
\begin{eqnarray*}
\frac{g(\bq(h))-g(\bq(0)) - h\rho_0(\bq)^\top M^{-1}p_0}{h^2} &=& \frac{1}2\left\{\ddot g(q_0) -\nabla^2 g(q_0)(M^{-1}p_0,M^{-1}p_0)+ O(h)\right\}, \\
\rho_0(\bq)^\top M^{-1}\rho_0(\bq) &=& \nabla g(q_0)^\top M^{-1} \nabla g(q_0) + O(h),\\
\rho_0(\bq)^\top M^{-1}\psi_0(\bq) &=& \nabla g(q_0)^\top M^{-1}\nabla U(q_0) + O(h),
\end{eqnarray*}
and this means that (\ref{consist}) tends to (\ref{lam0}), for $h\rightarrow 0$. Consequently, $\blam$ exists and is unique for all sufficiently small stepsizes $h>0$.\\

On the other hand, $g(q_1)-g(q_0) = g(q_1)=0$, provided that (see (\ref{hGamma}))
$$\hat\Gamma(\bq,\bp,\blam,h) = 0.$$
This means,
\begin{eqnarray}\label{sisl}
\lefteqn{\rho_0(\bq)^\top  M^{-1}p_0} \\ \nonumber
&=& h\sum_{j\ge0}\rho_j(\bq)^\top  M^{-1}\left\{
\xi_j[\psi_{j-1+\delta_{j0}}(\bq)+\rho_{j-1+\delta_{j0}}(\bq)\blam] - \xi_{j+1}[\psi_{j+1}(\bq)+\rho_{j+1}(\bq)\blam]\right\}\\ \nonumber
%%&&~+~h\xi_{s-1}\rho_{s-1}(\bq)^\top M^{-1}[\psi_{s-2}(\bq)+\rho_{s-2}(\bq)\blam]\\ \nonumber
&=& h\left\{ \left(\xi_0 \rho_0(\bq)^\top M^{-1}\rho_0(\bq) +\sum_{j\ge1} \xi_j\left[ \rho_j(\bq)^\top M^{-1}\rho_{j-1}(\bq) - \rho_{j-1}(\bq)^\top M^{-1}\rho_j(\bq)\right]\right)\blam \right.\\ \nonumber
&&~\left.+~\xi_0\rho_0(\bq)^\top M^{-1}\psi_0(\bq) +\sum_{j\ge1} \xi_j\left[ \rho_j(\bq)^\top M^{-1}\psi_{j-1}(\bq) - \rho_{j-1}(\bq)^\top M^{-1}\psi_j(\bq)\right] \right\},
\end{eqnarray}
and can be formally recast into the following linear system:
\begin{equation}\label{hMleqb}
A(h)\blam = b(h).
\end{equation}
Due to (\ref{cgj}) and Corollary~\ref{cor1}, the coefficient matrix reads: 
\begin{eqnarray}\label{Ah}
A(h) &=& h\xi_0 \rho_0(\bq)^\top M^{-1}\rho_0(\bq) +O(h^2) ~\equiv~ \frac{h}2 \nabla g(q_0)^\top M^{-1}\nabla g(q_0) + O(h^2),
\end{eqnarray}
and the right-hand side is
\begin{eqnarray}\nonumber
\lefteqn{b(h) ~=~ \rho_0(\bq)^\top  M^{-1}p_0 - \xi_0 h \rho_0(\bq)^\top M^{-1}\psi_0(\bq) + O(h^2)}  \\   \label{b}
&\equiv& \nabla g(q_0)^\top  M^{-1}p_0 + \frac{h}2\left[ \nabla^2 g(q_0) (M^{-1} p_0,M^{-1}p_0)- \nabla g(q_0)^\top M^{-1}\nabla U(q_0)\right] + O(h^2). \qquad
\end{eqnarray}
Consequently, (\ref{hMleqb}) is consistent with (\ref{lambda1}), since \,$ \nabla g(q_0)^\top  M^{-1}p_0=0$, thus giving
\begin{eqnarray*}
\lefteqn{\lambda_0 ~=~}\\
 &=& \left[\nabla g(q_0)^\top M^{-1}\nabla g(q_0) + O(h)\right]^{-1}\left[ \nabla^2 g(q_0) (M^{-1} p_0,M^{-1}p_0)- \nabla g(q_0)^\top M^{-1}\nabla U(q_0) + O(h)\right]\\
&\equiv& \lambda(q_0,p_0) + O(h).\end{eqnarray*}
From Theorem~\ref{hH_cons}, the conservation of energy follows. \\

The third statement of the theorem can be shown using the nonlinear variation of constants formula, by noting that for $t\in[0,h]$, $$\lambda(t)-\blam~\equiv~\underbrace{\lambda(q(t),p(t))-\lambda(q(0),p(0))}_{=O(h)} +O(h) ~=~O(h).$$ The last result follows from 
\begin{eqnarray*}
\nabla g(q_1)^\top M^{-1} p_1 &=& \nabla g(q(h)+O(h^2))^\top M^{-1} (p(h)+O(h^2)) \\[2mm]
&=& \underbrace{\nabla g(q(h))^\top M^{-1} p(h)}_{=0} + O(h^2) ~=~ O(h^2).\,\QED
\end{eqnarray*}

\bigskip
Next, let us consider the mesh
\begin{equation}\label{mesh}
t_n = nh, \quad n=0,\dots,N,   
\end{equation}
and the sequence of problems
\begin{equation}\label{constrp_lcn}
\dot \bq =  M^{-1}\bp, \qquad \dot \bp = -\nabla U(\bq)-\nabla g(\bq)\lambda_n, \qquad t\in[t_n,t_{n+1}],
\end{equation}
subject to initial conditions
\begin{equation}\label{unvn}
\bq(t_n) = q_n, \qquad \bp(t_n) = p_n,
\end{equation}
where $\lambda_n$ is a suitable constant vector. Then, the following result follows. 

\begin{theo}\label{induct}
Consider the IVPs (\ref{constrp_lcn})--(\ref{unvn}) and let us denote by $(q(t),p(t))$  the solution of the problem (\ref{constrp})--(\ref{q0p0}). Moreover, let us
assume that $(q_n,p_n)$ satisfies the following conditions: 
\begin{itemize}
\item $q_n-q(t_n) = O(h)$, \qquad $p_n-p(t_n)=O(h)$;
\item $\nabla g(q_n)^\top  M^{-1} \nabla g(q_n) \in\RR^{\nu\times \nu}$ is nonsingular;
\item $g(q_n) = 0$;
\item $\nabla g(q_n)^\top M^{-1} p_n = O(h^2)$.
\end{itemize}
Then, for all sufficiently small $h>0$,  $\exists!\lambda_n\in\RR^\nu$ such that the approximations
\begin{equation}\label{qn1pn1def}
q_{n+1}:=u(t_{n+1}), \qquad p_{n+1}:=v(t_{n+1}),
\end{equation}
satisfy:
\begin{itemize}
\item $g(q_{n+1})=0$ and, therefore, $H(q_{n+1},p_{n+1})=H(q_n,p_n)$;
\item $\lambda_n = \lambda(q(t_n),p(t_n)) + O(h)$;
\item $q_{n+1}-q(t_{n+1}) = O(h)$, \quad $p_{n+1}-p(t_{n+1})=O(h)$;
\item $\nabla g(q_{n+1})^\top M^{-1} p_{n+1} = O(h^2)$.
\end{itemize}
\end{theo}
\proof To show the first statement, we argue as we did to prove the first results in Theorem~\ref{lambda}. This yields $g(q_{n+1})=0$ provided that, cf.\ (\ref{hMleqb})--(\ref{b}),
$$A(h)\lambda_n = b(h)$$
where $A(h)$ and $b(h)$ are defined as in (\ref{Ah})--(\ref{b}) but with $q_0$ and $p_0$ replaced by $q_n$ and $p_n$, respectively. Consequently, from (\ref{lqp}), we obtain
\begin{eqnarray*}
\lefteqn{\lambda_n ~=~ \left[\nabla g(q_n)^\top M^{-1}\nabla g(q_n) + O(h)\right]^{-1}}\\
&&\left[ \nabla^2 g(q_n) (M^{-1} p_n,M^{-1}p_n)- \nabla g(q_n)^\top M^{-1}\nabla U(q_n)  -\frac{2}h\overbrace{\nabla g(q_n)^\top M^{-1} p_n}^{=O(h^2)} + O(h)\right]\\
&=& \left[\nabla g(q_n)^\top M^{-1}\nabla g(q_n) + O(h)\right]^{-1}\left[ \nabla^2 g(q_n) (M^{-1} p_n,M^{-1}p_n)- \nabla g(q_n)^\top M^{-1}\nabla U(q_n) + O(h)\right]\\[2mm]
&\equiv& \lambda(q_n,p_n) + O(h) ~=~\lambda(q(t_n)+O(h),p(t_n)+O(h))+O(h) ~=~\lambda(q(t_n),p(t_n))+O(h).\end{eqnarray*}
Energy conservation follows, as before, from Theorem~\ref{hH_cons}. Moreover, the nonlinear variation of constants formula, yields 
$$\pmatrix{c} q_{n+1}-q(t_{n+1})\\ p_{n+1}-p(t_{n+1}) \endpmatrix  = O(h) + O(h^2) = O(h),$$
due to $\lambda_n-\lambda(q(t),p(t))=O(h)$, for $t\in[t_n,t_{n+1}]$, and the hypothesis $q_n = q(t_n)+O(h)$, $p_n=p(t_n)+O(h)$.\\

In order to prove $\nabla g(q_{n+1})^\top M^{-1} p_{n+1} = O(h^2)$, we note that from the hypothesis
$$\nabla g(q_n)^\top M^{-1}p_n = O(h^2),$$ 
the existence of $\tp_n\in\RR^m$ such that
$$p_n-\tp_n=O(h^2), \qquad \nabla g(q_n)^\top M^{-1} \tp_n = 0$$
follows. Using $(q_n,\tp_n)$ as local initial conditions for (\ref{constrp_lcn}), and repeating above steps to satisfy the constraints at $t_{n+1}$, we obtain  
$$\tlam_n = \lambda(q_n,\tp_n) + O(h) \equiv \lambda(q_n,p_n+O(h^2)) + O(h) = \lambda(q_n,p_n) + O(h) \equiv \lambda_n+O(h),$$
and corresponding approximations $\tq_{n+1}$, $\tp_{n+1}$ such that
$$g(\tq_{n+1})=0, \qquad \nabla g(\tq_{n+1})^\top M^{-1} \tp_{n+1} = O(h^2).$$
From $p_n-\tp_n = O(h^2)$ and $\lambda_n-\tlam_n=O(h)$, the nonlinear variation of constants formula yields, 
$$q_{n+1}-\tq_{n+1} = O(h^2), \qquad p_{n+1}-\tp_{n+1}=O(h^2).$$
Consequently, 
\begin{eqnarray*}
\nabla g(q_{n+1})^\top M^{-1} p_{n+1} &=& \nabla g(\tq_{n+1}+O(h^2))^\top M^{-1} (\tp_{n+1}+O(h^2)) \\[2mm]
&=& \underbrace{\nabla g(\tq_{n+1})^\top M^{-1} \tp_{n+1}}_{=O(h^2)}~+~O(h^2) ~=~ O(h^2).\,\QED
\end{eqnarray*}

\bigskip
By means of Theorem~\ref{induct}, a straightforward induction argument enables to show the following relaxed version of Theorem~\ref{capdcazz}.

\begin{cor}\label{capdcazz1}
For all sufficiently small stepsizes $h>0$, the above procedure defines a sequence of approximations $(q_n,p_n)$ such that,
for all $n=1,2,\ldots,$
\begin{equation}\label{ordqp1}
q_n = q(nh) + O(h), \qquad p_n = p(nh) + O(h), \qquad g(q_n)=0, \qquad \nabla g(q_n)^\top M^{-1}p_n = O(h^2).
\end{equation}
Moreover, $(q_{n+1},p_{n+1})$ is obtained from $(q_n,p_n)$ utilizing a constant vector $\lambda_n$ such that
$$%\begin{equation}\label{ordl}
\lambda_n = \lambda( q(nh), p(nh) ) + O(h),
$$%\end{equation}
where $\lambda(q,p)$ is the function defined in (\ref{lqp}) and consequently, ~$H(q_{n+1},p_{n+1}) = H(q_n,p_n)$ holds.
\end{cor}

\medskip
The fact that (\ref{ordqp}) holds, in place of the weaker result (\ref{ordqp1}), is due to the {\em symmetry} of the proposed procedure. Note that a symmetric method is necessarily of even order \cite[Theorem\,3.2]{HLW2006}. The method is symmetric since we have shown that, for  all sufficiently small $h>0$, there exists a unique $\lambda_n$ such that from the solution of (\ref{constrp_lcn})--(\ref{unvn}), with
\begin{equation}\label{tn}
(q_n,p_n), \quad \quad g(q_n)=0, %\qquad \nabla g(q_n)^\top M^{-1} p_n = O(h^2),
\end{equation}
we arrive at a new point, where 
\begin{equation}\label{tn1}
(q_{n+1},p_{n+1}), \quad  \quad g(q_{n+1})=0, %\qquad \nabla g(q_{n+1})^\top M^{-1} p_{n+1} = O(h^2).
\end{equation}
Since $\lambda_n$ is unique, when we start from (\ref{tn1}) and solve backward in time (\ref{constrp_lcn}), we arrive at (\ref{tn}), i.e.\ the procedure is symmetric. As a consequence of the symmetry of the method, the approximation order of $(q_n,p_n)$ is even and therefore, (\ref{ordqp}) holds in place of (\ref{ordqp1}). This completes the proof of Theorem~\ref{capdcazz}.

\bigskip
In the next theorem, we summarize in a more comprehensive form the statements derived previously       in this section.

\begin{theo}\label{capdcazz2}
Let us consider the problem (\ref{constrp})--(\ref{consist0}), the mesh (\ref{mesh}), and the sequence of problems (\ref{constrp_lcn})--(\ref{unvn}). The constant vector $\lambda_n$ is chosen in such a way that for the new approximations defined by (\ref{qn1pn1def}), $g(q_{n+1})=0$ follows. Then, for all  sufficiently small stepsizes $h>0$, the above procedure defines a sequence of approximations
$(q_n,p_n,\lambda_n)$ satisfying\,\footnote{For the definition of $\lambda(q,p)$ see (\ref{lqp}).}
\begin{eqnarray*}
q_n &=& q(nh) + O(h^2), \\
p_n &=& p(nh) + O(h^2), \\
g(q_n) &=& 0, \\
\nabla g(q_n)^\top M^{-1}p_n &=& O(h^2),\\
\lambda_n &=& \lambda( q(nh), p(nh) ) + O(h),\\
H(q_n,p_n) &=& H(q_0,p_0), \qquad n=0,1,\dots,N.
\end{eqnarray*}
Moreover, in case that $\lambda$ is constant, $\lambda(q(t),p(t))\equiv \bar\lambda$, $\forall t\in[0,T]$, the following statements hold: 
\begin{eqnarray*}
q_n &=& q(nh), \\
p_n &=& p(nh), \\
g(q_n) &=& 0, \\
\nabla g(q_n)^\top M^{-1}p_n &=& 0,\\
\lambda_n &=& \bar\lambda,\\
H(q_n,p_n) &=& H(q_0,p_0), \qquad n=0,1,\dots,N.
\end{eqnarray*}
\end{theo}

\bigskip
With other words, the discrete solution is {\em exact}, when the vector of the Lagrange multipliers is constant. Otherwise, it is second-order accurate for $(q_n,p_n)$ and first order accurate for $\lambda_n$. In the latter case, the constraints and the Hamiltonian are conserved, while the hidden constraints remain $O(h^2)$ close to zero.

\section{Polynomial approximation}\label{poly}

The first step towards the numerical solution of (\ref{constrp})--(\ref{consist0}) is to truncate the series in the right-hand side of (\ref{constrp1}),
\begin{equation}\label{constrps}
\dot \bq(ch) = \sum_{j=0}^{s-1} P_j(c) \gamma_j(\bp), \qquad \dot \bp(ch) = -\sum_{j=0}^{s-1} P_j(c)[ \psi_j(\bq) +\rho_j(\bq)\blam], \qquad  c\in[0,1].
\end{equation}
Here, the coefficients  $\gamma_j, \psi_j$, and $\rho_j$ are as those defined in (\ref{cgj}). By imposing the initial conditions, the local approximation over the first step becomes
\begin{equation}\label{qpchs}
\bq(ch) = q_0 + h\sum_{j=0}^{s-1} \int_0^cP_j(x)\dd x\, \gamma_j(\bp), \quad  \bp(ch) = p_0 -h\sum_{j=0}^{s-1} \int_0^cP_j(x)\dd x[ \psi_j(\bq) +\rho_j(\bq)\blam], \quad  c\in[0,1],
\end{equation}
with the new approximations which are formally still given by (\ref{q1p1def}).
Again, the constant vector of the multipliers is uniquely determined by requiring that the constraints are satisfied at $t_1=h$. Consequently, we have the following expression, in place of (\ref{sisl}):
\begin{eqnarray}\label{sisls}
\lefteqn{\rho_0(\bq)^\top  M^{-1}p_0} \\ \nonumber
&=& h\sum_{j=0}^{s-2}\rho_j(\bq)^\top  M^{-1}\left\{
\xi_j[\psi_{j-1+\delta_{j0}}(\bq)+\rho_{j-1+\delta_{j0}}(\bq)\blam] - \xi_{j+1}[\psi_{j+1}(\bq)+\rho_{j+1}(\bq)\blam]\right\}\\ \nonumber
&&~+~h\xi_{s-1}\rho_{s-1}(\bq)^\top M^{-1}[\psi_{s-2}(\bq)+\rho_{s-2}(\bq)\blam]\\ \nonumber
&=& h\left\{ \left(\xi_0 \rho_0(\bq)^\top M^{-1}\rho_0(\bq) +\sum_{j=1}^{s-1} \xi_j\left[ \rho_j(\bq)^\top M^{-1}\rho_{j-1}(\bq) - \rho_{j-1}(\bq)^\top M^{-1}\rho_j(\bq)\right]\right)\blam \right.\\ \nonumber
&&~\left.+~\xi_0\rho_0(\bq)^\top M^{-1}\psi_0(\bq) +\sum_{j=1}^{s-1} \xi_j\left[ \rho_j(\bq)^\top M^{-1}\psi_{j-1}(\bq) - \rho_{j-1}(\bq)^\top M^{-1}\psi_j(\bq)\right] \right\},
\end{eqnarray}
which, in turn, yields equations which are formally similar to (\ref{hMleqb})--(\ref{b})\footnote{As before, this basic step defines a symmetric procedure.}.
The process is then repeated by defining the mesh (\ref{mesh}) and considering the local problems
\begin{eqnarray}\nonumber
\dot \bq_n(ch) = \sum_{j=0}^{s-1} P_j(c) \gamma_j(\bp_n), && \dot \bp_n(ch) ~=~ -\sum_{j=0}^{s-1} P_j(c)[ \psi_j(\bq_n) +\rho_j(\bq_n)\lambda_n], \quad  c\in[0,1],\\
\label{constrp_lcns}
u_n(0)=q_n, && v_n(0) = p_n, \quad n=0,\dots,N-1, 
\end{eqnarray}
where the coefficients $\gamma_j(\bp_n), \psi_j(\bq_n), \rho_j(\bq_n)$ are defined in (\ref{cgj}), with $\bq$ and $\bp$ replaced by $\bq_n$ and $\bp_n$, respectively.
Consequently, we formally obtain the piecewise polynomial approximation,
\begin{eqnarray}\label{ppapp}
\bq_n(ch) &=& q_n+h\sum_{j=0}^{s-1} \int_0^cP_j(x)\dd x \gamma_j(\bp_n), \\ \nonumber
\bp_n(ch)  &=& p_n -h\sum_{j=0}^{s-1} \int_0^cP_j(x)\dd x[ \psi_j(\bq_n) +\rho_j(\bq_n)\lambda_n], \qquad  c\in[0,1],
\end{eqnarray}
with the new approximations given by (see (\ref{orto}))
\begin{equation}\label{qn1pn1defs}
q_{n+1} := \bq_n(h) \equiv q_n+h\gamma_0(\bp_n), \qquad
p_{n+1} := \bp_n(h) \equiv p_n-h[\psi_0(\bq_n) +\rho_0(\bq_n)\lambda_n].
\end{equation}
As before, the constant vector $\lambda_n\in\RR^\nu$ is chosen to satisfy the constraints $g(q_{n+1})=0$ and it is implicitly defined by the equation,
\begin{eqnarray}\label{sislsn}
\lefteqn{\rho_0(\bq_n)^\top  M^{-1}p_n} \\ \nonumber
&=& h\left\{ \left(\xi_0 \rho_0(\bq_n)^\top M^{-1}\rho_0(\bq_n) +\sum_{j=1}^{s-1} \xi_j\left[ \rho_j(\bq_n)^\top M^{-1}\rho_{j-1}(\bq_n) - \rho_{j-1}(\bq_n)^\top M^{-1}\rho_j(\bq_n)\right]\right)\lambda_n \right.\\ \nonumber
&&~\left.+~\xi_0\rho_0(\bq_n)^\top M^{-1}\psi_0(\bq_n) +\sum_{j=1}^{s-1} \xi_j\left[ \rho_j(\bq_n)^\top M^{-1}\psi_{j-1}(\bq_n) - \rho_{j-1}(\bq_n)^\top M^{-1}\psi_j(\bq_n)\right] \right\}.
\end{eqnarray}
This equation reduces to (\ref{sisls}) for $n=0$. Using arguments similar to those from the previous section (see also \cite{BIT2012_2}), it is possible to show the following result. This result is a counterpart to Theorem~\ref{capdcazz2} for the piecewise polynomial approximation (\ref{ppapp}) to the solution $(q(t),p(t))$ of problem (\ref{constrp})--(\ref{consist0}).

\begin{theo}\label{capdcazz3}
For all sufficiently small stepsizes $h>0$, the approximation procedure (\ref{constrp_lcns})--(\ref{sislsn}) is well defined and provides a sequence of approximations  $(q_n,p_n,\lambda_n)$ such that
\begin{eqnarray*}
q_n &=& q(nh) + O(h^2), \\
p_n &=& p(nh) + O(h^2), \\
g(q_n) &=& 0, \\
\nabla g(q_n)^\top M^{-1}p_n &=& O(h^2),\\
\lambda_n &=& \lambda( q(nh), p(nh) ) + O(h),\\
H(q_n,p_n) &=& H(q_0,p_0), \qquad n=0,1,\dots,N.
\end{eqnarray*}
Moreover, in case that $\lambda$ is constant, $\lambda(q(t),p(t))\equiv \bar\lambda$, $\forall t\in[0,T]$, the following statements hold: 
\begin{eqnarray*}
q_n &=& q(nh) + O(h^{2s}), \\
p_n &=& p(nh) + O(h^{2s}), \\
g(q_n) &=& 0, \\
\nabla g(q_n)^\top M^{-1}p_n &=& O(h^{2s}),\\
\lambda_n &=& \bar\lambda +O(h^{2s}),\\
H(q_n,p_n) &=& H(q_0,p_0), \qquad n=0,1,\dots,N.
\end{eqnarray*}
\end{theo}

\section{Full discretization}\label{HBVMs}

In order to cast the above algorithm into a computational method, the integrals defining the coefficients $\gamma_j(\bp), \, \psi_j(\bq), \, \rho_j(\bq), \, j=0,\dots,s-1,$ in (\ref{constrps}), need to be approximated.\footnote{Since the method is a one-step method, we shall, as usual, only consider the first step.} To this aim, following the discussion in \cite{BIT2010,BIT2012_2,LIMbook2016}, we use the Gauss-Legendre quadrature of order $2k$ (the interpolatory quadrature formula based at the zeros of $P_k(c)$), with nodes and weights $(\hc_i,\hb_i)$, where $k\ge s$.
Consequently, 
\begin{eqnarray}\nonumber
\gamma_j(\bp) &\approx& \hat\gamma_j ~:=~  M^{-1}\sum_{\ell=1}^k \hb_\ell P_j(\hc_\ell)\bp(\hc_\ell h), \qquad \psi_j(\bq) ~\approx~ \hat\psi_j ~:=~ \sum_{\ell=1}^k \hb_\ell  P_j(\hc_\ell)\nabla U(\bq(\hc_\ell h)),\\  \label{cgjd}
\rho_j(\bq) &\approx& \hat\rho_j ~:=~ \sum_{\ell=1}^k \hb_\ell  P_j(\hc_\ell)\nabla g(\bq(\hc_\ell h)), \qquad j=0,\dots,s-1.
\end{eqnarray}
Formally, this is a $k$-stage Runge-Kutta method, whose computational cost depends on $s$ rather than on $k$, since the actual unknowns are the $3s$ coefficients (\ref{cgjd}) and the vector $\blam$ (see (\ref{constrps})). We refer, to \cite{BIT2011,LIMbook2016} for details.\\

Let us now formulate the discrete problem to be solved in each integration step. We first define the matrices, cf.\ (\ref{PsIsOm}),
$$\hat\P_s = \left( P_{j-1}(\hc_i)\right), ~\, \hat\I_s = \left( \int_0^{\hc_i} P_{j-1}(x)\dd x\right)~\in \RR^{k\times s}, \qquad \hat\Omega = \diag(\hb_1,\dots,\hb_k)\in\RR^{k\times k},$$
and the vectors and matrices
$$e = \pmatrix{c} 1\\ \vdots\\ 1\endpmatrix \in\RR^k, \qquad \hat\bfgamma = \pmatrix{c} \hat\gamma_0\\ \vdots \\ \hat\gamma_{s-1}\endpmatrix,~\, \hat\bfpsi = \pmatrix{c} \hat\psi_0\\ \vdots \\ \hat\psi_{s-1}\endpmatrix ~\in\RR^{sm}, \qquad \hat\bfrho = \pmatrix{c} \hat\rho_0\\ \vdots \\ \hat\rho_{s-1}\endpmatrix\in\RR^{sm\times \nu}.$$
Recall that $m$ is the dimension of the continuous problem and $\nu$ is the number of constraints.\\

Then, the $3s$ equations from (\ref{cgjd}), defining the discrete problem to be solved, amount to the system of equations, of (block) dimension $s$,
\begin{eqnarray}\nonumber
\hat\bfgamma &=& \hat\P_s^\top\hat\Omega \otimes  M^{-1} \left[ e\otimes p_0 - h\hat\I_s\otimes I_m\left( \hat\bfpsi +\hat\bfrho\blam\right) \right],\\ \label{dispro}
\hat\bfpsi &=& \hat\P_s^\top\hat\Omega \otimes I_m \nabla U\left( e\otimes q_0 +h\hat\I_s\otimes I_m \hat\bfgamma\right),\\ \nonumber
\hat\bfrho &=& \hat\P_s^\top\hat\Omega \otimes I_m \nabla g\left( e\otimes q_0 +h\hat\I_s\otimes I_m \hat\bfgamma\right).
\end{eqnarray}
We augment (\ref{dispro}) by the equation (\ref{sislsn}) for $\blam$ which, by taking (\ref{cgjd}) into account, can be rewritten as
\begin{eqnarray}\nonumber
\lefteqn{h\left[\xi_0 \hat\rho_0^\top M^{-1}\hat\rho_0 +\sum_{j=1}^{s-1} \xi_j\left( \hat\rho_j^\top M^{-1}\hat\rho_{j-1} - \hat\rho_{j-1}^\top M^{-1}\hat\rho_j\right)\right]\blam}\\  \label{sisld}
&=&\hat\rho_0^\top  M^{-1} \left(p_0-h\xi_0\hat\psi_0\right)-h\sum_{j=1}^{s-1} \xi_j\left( \hat\rho_j^\top M^{-1}\hat\psi_{j-1} - \hat\rho_{j-1}^\top M^{-1}\hat\psi_j\right).
\end{eqnarray}
In (\ref{dispro}), $\nabla U$, when evaluated in a block vector of (block) dimension $k$, stands for the block vector made up of the $k$ vectors resulting from the corresponding application of the function. The same straightforward notation is used for $\nabla g$.
The new approximation is then given by, see (\ref{q1p1def}) and (\ref{cgjd}),
\begin{equation}\label{y1} q_1 = q_0 + h\hat\gamma_0, \qquad p_1 = p_0 -h[ \hat\psi_0 +\hat\rho_0\blam].\end{equation}
Note that the equations in (\ref{dispro}), together with (\ref{y1}), formally {\em coincide} with those provided by a HBVM$(k,s)$ method\footnote{Here, $s$ is the degree of the polynomial approximation and $k$ defines the order (actually equal to $2k$) of the quadrature in the approximations (\ref{cgjd}).} applied to solve the problem defined by the Hamiltonian (\ref{Hqpl}), where the vector of the multiplier $\blam$ is considered as a parameter,%\footnote{Clearly, the problem (\ref{nogq}) represents the ``unconstrained'' version of (\ref{constrp}).}
$$%\begin{equation}\label{nogq}
\dot q =  M^{-1}p, \quad \dot p = -\nabla U(q)-\nabla g(q)\blam, \quad  t\ge0, \qquad q(0)=q_0, ~p(0)=p_0,
$$%\end{equation}
cf.\ \cite{BIT2011,LIMbook2016} for details. Consequently, equation (\ref{sisld}) defines the proper extension for handling the constrained Hamiltonian problem (\ref{Hqp})--(\ref{gq}). For this reason, we continue to refer to the numerical method specified in (\ref{dispro})--(\ref{y1}) as to HBVM$(k,s)$. Now, it is a {\em ready to use} numerical procedure. \\

The discrete problem (\ref{dispro})--(\ref{sisld}) can be solved via a straightforward fixed-point iteration, which converges under regularity assumptions, for all sufficiently small stepsizes $h>0$.\footnote{We refer to \cite{BIT2011,BFCI2014,LIMbook2016} for further details on procedures for solving the involved discrete problems. They are based on suitable Newton-splitting procedures, already implemented in computational codes \cite{BrMa2002,BrMa2004,BrMa2007}.} Moreover, for separable Hamiltonians, as it is the case in (\ref{Hqp}), the last two equations in (\ref{dispro}) can be substituted into the first one, resulting in a {\em single} vector equation for $\hat\bfgamma$. By setting
\begin{equation}\label{teta}
\Theta_\blam(q) := \nabla U (q)+\nabla g(q)\blam,
\end{equation}
we obtain
\begin{equation}\label{sologamma}
\hat\bfgamma = \hat\P_s^\top\hat\Omega \otimes  M^{-1} \left[ e\otimes p_0 - h\hat\I_s\hat\P_s^\top\hat\Omega\otimes I_m \Theta_\blam\left( e\otimes q_0 +h\hat\I_s\otimes I_m \hat\bfgamma\right) \right],
\end{equation}
plus (\ref{sisld}) for $\blam$.\footnote{Clearly, $\hat\bfpsi$ and $\hat\bfrho$ can be computed via the last two equations in (\ref{dispro}), once $\hat\bfgamma$ and $\blam$ are known.}  We skip further details, since they are exactly the same as for the original HBVMs, when applied to solve separable (unconstrained) Hamiltonian problems \cite{BIT2011,LIMbook2016}.\\

The following result follows from Theorem~\ref{capdcazz3} along with standard arguments from the analysis of HBVMs \cite{BIT2012_2,LIMbook2016}.\footnote{For brevity, we omit the proof here.}

\begin{theo}\label{capdcazz4}
For all sufficiently small stepsizes $h>0$, the HBVM$(k,s)$ method (\ref{sisld})--(\ref{sologamma}) is well defined and symmetric. It provides a sequence of approximations  $(q_n,p_n,\lambda_n)$, $n=0,1,\dots,N$, such that 
\begin{eqnarray*}
q_n &=& q(nh) + O(h^2), \\
p_n &=& p(nh) + O(h^2), \\
\nabla g(q_n)^\top M^{-1}p_n &=& O(h^2),\\
\lambda_n &=& \lambda( q(nh), p(nh) ) + O(h),
\end{eqnarray*}
and
\begin{eqnarray}\nonumber
g(q_n) &=&  \left\{ \begin{array}{cl} 0, &\mbox{if~$g$~is~a~polynomial~of~degree~not~larger~than~$2k/s$,}\\[2mm]
O(h^{2k}), &\mbox{otherwise},\end{array}\right.\\[-2mm]
\label{Hg}\\ \nonumber
H(q_n,p_n)-H(q_0,p_0) &=&
\left\{ \begin{array}{cl} 0, &\mbox{if~$H$~is~a~polynomial~of~degree~not~larger~than~$2k/s$,}\\[2mm]
O(h^{2k}), &\mbox{otherwise}.\end{array}\right.
\end{eqnarray}
Moreover, in case $\lambda$ is constant, $\lambda(q(t),p(t))\equiv \bar\lambda$, $\forall t\in[0,T]$, the following statements hold:
\begin{eqnarray*}
q_n &=& q(nh) + O(h^{2s}), \\
p_n &=& p(nh) + O(h^{2s}), \\
%g(q_n) &=& 0, \\
\nabla g(q_n)^\top M^{-1}p_n &=& O(h^{2s}),\\
\lambda_n &=& \bar\lambda +O(h^{2s}).
%H(q_n,p_n) &=& H(q_0,p_0), \qquad n=0,1,\dots,N.
\end{eqnarray*}
\end{theo}

\begin{rem} We stress that by (\ref{Hg}), an exact or a (at least) practical conservation of \underline{both} the constraints and the Hamiltonian can \underline{always} be guaranteed. In fact, by choosing sufficiently large $k$, either the quadrature becomes    exact, in the polynomial case, or the quadrature error is within the round-off error level, in the non polynomial case. This feature of the method will be always exploited in the numerical tests discussed in Section~\ref{numtest}.
\end{rem}

Finally, we shall mention that for $k=s$, the HBVM$(s,s)$ method reduces to the $s$-stage Gauss collocation method, \cite{BIT2010,BIT2012_2,LIMbook2016}, which is symplectic. Moreover, in the limit $k\rightarrow\infty$, we retrieve the formulae studied in Section~\ref{poly}. This means that our approach can be also considered in the framework of Runge-Kutta methods with continuous stages \cite{BIT2010,H2010}.

\section{Numerical tests}\label{numtest}
In this section, to illustrate the theoretical properties of HBVM$(k,s)$, we apply them to numerically simulate some Hamiltonian problems of the form (\ref{Hqp})--(\ref{gq}) with holonomic constraints. In the focus of our attention are properties described in Theorem~\ref{capdcazz4}.  % code: hbvm_hc.m

\subsection{Pendulum} We begin with the {\em planar pendulum} in Cartesian coordinates, where a massless rod of length $L$ connects a point of mass $m$ to a fixed point (the origin). We assume a unit mass and length, $m=1$ and $L=1$, and normalize the gravity acceleration. Then, the Hamiltonian is given by
\begin{equation}\label{pend_H}
H(q,p) = \frac{1}2 p^\top p + e_2^\top q, \qquad e_2 := \pmatrix{c} 0\\ 1\endpmatrix, \quad q := \pmatrix{c} x\\ y\endpmatrix, \quad p:= \dot q\,\in\RR^2,
\end{equation}
and is subject to the constraint
\begin{equation}\label{pend_c}
g(q) \equiv q^\top q -1 = 0.
\end{equation}
We also prescribe the initial conditions of the form
\begin{equation}\label{pend_y0}
q(0) = (0,\,-1)^\top, \qquad p(0) = (1,\,0)^\top.
\end{equation}
Consequently, the constrained Hamiltonian problem reads: 
\begin{eqnarray}\label{pend_c1}
&&\ddot x =  -2x\lambda, \qquad \ddot y = -1-2y\lambda, \qquad x^2+y^2=1, \\[2mm] \nonumber
&& x(0)=0,\quad y(0)=-1,\quad \dot x(0) = 1,\quad \dot y(0) = 0.
\end{eqnarray}
In order to obtain a reference solution, we rewrite the problem in polar coordinates in such a way that $\theta=0$ locates the pendulum at its stable rest position, so that 
\begin{equation}\label{xy}
x = \sin\theta, \qquad y = -\cos\theta.
\end{equation}
Thus, we arrive at the unconstrained Hamiltonian problem
\begin{equation}\label{pend_u}
\ddot\theta +\sin\theta = 0, \qquad \theta(0) = 0, \quad \dot\theta(0) = 1.
\end{equation}
Once this problem is solved, the solution of (\ref{pend_c1}) is recovered via the transformations (\ref{xy}). Moreover, the Lagrange multiplier in (\ref{pend_c}) turns out to be given by
\begin{equation}\label{pend_l}
\lambda = \frac{1}2\left( \dot\theta^2 +\cos\theta\right).
\end{equation}
To compute the reference solution for (\ref{pend_c1}), we solve (\ref{pend_u}) by means of a HBVM$(12,6)$ method\,\footnote{For unconstrained Hamiltonian problems.} of order 12 which is practically energy-conserving. \\

According to (\ref{pend_H}) and (\ref{pend_c}), the Hamiltonian and the constraint are quadratic, so we expect HBVM$(s,s)$ to conserve the energy and the constraint. In Table~\ref{tab0}, we list the following quantities, obtained from the HBVM$(s,s)$ methods for $s=1,2,3$, the stepsizes $h=10^{-1}2^{-n}$ and the interval of integration $[0,10]$:
\begin{itemize}
\item the solution error ($e_s$),
\item the multiplier error ($e_\lambda$),
\item the Hamiltonian error $(e_H)$,
\item the constraint error $(e_g)$;
\item the {\em hidden constraint} error, defined by 
$$e_{hc} := \max_n\,  2|x_n\dot{x}_n+y_n\dot{y}_n|.$$
\end{itemize}
As predicted in Theorem~\ref{capdcazz4}, we can see that
\begin{itemize}
\item all methods are second-order accurate in the space variables, with HBVM(1,1) less accu\-ra\-te than the others;
\item all methods are first-order accurate in the Lagrange multiplier;
\item all methods exactly conserve the Hamiltonian and the constraint;
\item all methods are second-order accurate in the hidden constraint.
\end{itemize}

\begin{table}[t]
\caption{\label{tab0} Planar pendulum (\ref{pend_H})--(\ref{pend_c}). Errors from HBVM$(s,s)$ method for $s=1,2,3$, when solving the problem over the interval $[0,10]$ with stepsizes $h=10^{-1}2^{-n}$.}  % see erroconic.m
\smallskip
\centerline{\begin{tabular}{|r|rrrrrrrr|}
\hline
\hline
       \multicolumn{9}{|c|}{$s=1$}\\
       \hline
$n$ & $e_s$  & rate  & $e_\lambda$ & rate  & $e_H$ & $e_g$ & $e_{hc}$ & rate   \\
\hline
  0 & 2.5700e-02 & -- & 3.4253e-02 & -- & 5.5511e-17 & 9.9920e-16 & 2.3487e-03 & --\\
  1 & 6.4260e-03 & 2.00 & 1.7386e-02 & 0.98 & 1.1102e-16 & 6.6613e-16 & 5.8639e-04 & 2.00\\
  2 & 1.6070e-03 & 2.00 & 8.7406e-03 & 0.99 & 1.1102e-16 & 5.5511e-16 & 1.4654e-04 & 2.00\\
  3 & 4.0181e-04 & 2.00 & 4.3835e-03 & 1.00 & 1.1102e-16 & 1.9984e-15 & 3.6633e-05 & 2.00\\
  4 & 1.0045e-04 & 2.00 & 2.1948e-03 & 1.00 & 1.1102e-16 & 1.8874e-15 & 9.1580e-06 & 2.00\\
  5 & 2.5114e-05 & 2.00 & 1.0982e-03 & 1.00 & 1.1102e-16 & 2.4425e-15 & 2.2895e-06 & 2.00\\
  6 & 6.2785e-06 & 2.00 & 5.4929e-04 & 1.00 & 1.1102e-16 & 3.4417e-15 & 5.7238e-07 & 2.00\\
  7 & 1.5696e-06 & 2.00 & 2.7470e-04 & 1.00 & 1.1102e-16 & 6.1062e-15 & 1.4311e-07 & 2.00\\
  8 & 3.9248e-07 & 2.00 & 1.3743e-04 & 1.00 & 1.1102e-16 & 8.5487e-15 & 3.5902e-08 & 2.00\\
\hline
\hline
       \multicolumn{9}{|c|}{$s=2$}\\
       \hline
$n$ & $e_s$  & rate  & $e_\lambda$ & rate  & $e_H$ & $e_g$ & $e_{hc}$ & rate   \\
\hline
  0 & 1.6695e-03 & -- & 3.5176e-02 & -- & 1.1102e-16 & 8.8818e-16 & 2.3539e-03 & --\\
  1 & 4.1412e-04 & 2.01 & 1.7585e-02 & 1.00 & 1.1102e-16 & 8.8818e-16 & 5.8670e-04 & 2.00\\
  2 & 1.0332e-04 & 2.00 & 8.7919e-03 & 1.00 & 1.1102e-16 & 1.1102e-15 & 1.4656e-04 & 2.00\\
  3 & 2.5816e-05 & 2.00 & 4.3958e-03 & 1.00 & 1.1102e-16 & 9.9920e-16 & 3.6634e-05 & 2.00\\
  4 & 6.4533e-06 & 2.00 & 2.1979e-03 & 1.00 & 1.1102e-16 & 1.8874e-15 & 9.1581e-06 & 2.00\\
  5 & 1.6133e-06 & 2.00 & 1.0990e-03 & 1.00 & 1.1102e-16 & 2.7756e-15 & 2.2895e-06 & 2.00\\
  6 & 4.0332e-07 & 2.00 & 5.4948e-04 & 1.00 & 1.1102e-16 & 3.6637e-15 & 5.7238e-07 & 2.00\\
  7 & 1.0086e-07 & 2.00 & 2.7477e-04 & 1.00 & 1.1102e-16 & 5.5511e-15 & 1.4314e-07 & 2.00\\
  8 & 2.5286e-08 & 2.00 & 1.3751e-04 & 1.00 & 1.1102e-16 & 1.0547e-14 & 3.5884e-08 & 2.00\\
\hline
\hline
       \multicolumn{9}{|c|}{$s=3$}\\
       \hline
$n$ & $e_s$  & rate  & $e_\lambda$ & rate  & $e_H$ & $e_g$ & $e_{hc}$ & rate   \\
\hline
  0 & 1.6658e-03 & -- & 3.5178e-02 & -- & 1.1102e-16 & 1.1102e-15 & 2.3539e-03 & --\\
  1 & 4.1386e-04 & 2.01 & 1.7585e-02 & 1.00 & 1.1102e-16 & 8.8818e-16 & 5.8670e-04 & 2.00\\
  2 & 1.0331e-04 & 2.00 & 8.7919e-03 & 1.00 & 1.1102e-16 & 7.7716e-16 & 1.4656e-04 & 2.00\\
  3 & 2.5815e-05 & 2.00 & 4.3958e-03 & 1.00 & 1.1102e-16 & 8.8818e-16 & 3.6634e-05 & 2.00\\
  4 & 6.4532e-06 & 2.00 & 2.1979e-03 & 1.00 & 1.1102e-16 & 1.6653e-15 & 9.1581e-06 & 2.00\\
  5 & 1.6133e-06 & 2.00 & 1.0990e-03 & 1.00 & 1.1102e-16 & 2.9976e-15 & 2.2895e-06 & 2.00\\
  6 & 4.0331e-07 & 2.00 & 5.4948e-04 & 1.00 & 1.1102e-16 & 3.4417e-15 & 5.7238e-07 & 2.00\\
  7 & 1.0086e-07 & 2.00 & 2.7477e-04 & 1.00 & 1.1102e-16 & 5.2180e-15 & 1.4314e-07 & 2.00\\
  8 & 2.5220e-08 & 2.00 & 1.3739e-04 & 1.00 & 1.1102e-16 & 8.9928e-15 & 3.5791e-08 & 2.00\\
\hline
\hline
\end{tabular}}
\end{table}

\subsection{Conical pendulum}  % conicalpend.m
Next, we consider the so-called {\em conical pendulum}, a particular case of the {\em spherical pendulum}, namely a pendulum of mass $m$, which is connected to a fixed point (i.e., the origin) by a massless rod of length $L$. For the conical pendulum, the initial condition is chosen in such a way that the motion is periodic with period $T$ and occurs in the horizontal plane $q_3 =z_0$\footnote{Clearly, $0>z_0>-L$.}. Again, assuming $m=1$ and $L=1$, and normalizing the acceleration of gravity, the Hamiltonian is
\begin{equation}\label{cpend_H}
H(q,p) := \frac{1}2p^\top p + e_3^\top q, \qquad e_3 := \pmatrix{c} 0\\ 0\\ 1\endpmatrix, \quad  q := \pmatrix{c} x\\ y\\ z\endpmatrix, \quad p := \dot q\,\in\RR^3,
\end{equation}
with the constraint
\begin{equation}\label{cpend_c}
g(q) := q^\top q-1 = 0.
\end{equation}
Here, we prescribe the consistent initial conditions
\begin{equation}\label{cpend_y0}
q(0) = \pmatrix{r} 2^{-\frac{1}2}\\ 0\\ -2^{-\frac{1}2} \endpmatrix, \qquad
p(0) = \pmatrix{r} 0 \\ 2^{-\frac{1}4} \\ 0 \endpmatrix,
\end{equation}
generating a periodic motion with
\begin{equation}\label{Tz0}
T = 2^{\frac{3}4}\pi, \qquad z_0 = -2^{-\frac{1}2}.
\end{equation}
Moreover, in such a case, the multiplier $\blam$, which has the physical meaning of the tension on the rod, has to be constant and is given by 
\begin{equation}\label{lambda0}
\lambda_0 = 2^{-\frac{1}2}.
\end{equation}
Note that the Hamiltonian and the constraint are quadratic and according to Theorem~\ref{capdcazz4}, any HBVM$(s,s)$ method conserves both of them and has order $2s$.
In Table~\ref{tab1}, we list the errors in
\begin{itemize}
\item the solution ($e_s$),
\item the multiplier ($e_\lambda$),
\item the Hamiltonian $(e_H)$,
\item the constraints $(e_g)$;
\item the {\em hidden constraints}, defined by 
\begin{equation}\label{ehc}
e_{hc} := \max_n \| \nabla g(q_n)^\top M^{-1}p_n\|.
\end{equation}
\end{itemize}
The problem is solved over $10$ periods, with stepsizes $h=T/n$. As expected, the estimated rate of convergence for HBVM$(s,s)$, $s=1,2,3,4$, is  $2s$. Also, the Hamiltonian and the constraint are conserved up to round-off errors.  Remarkably, also the error in the multiplier $(e_\lambda)$ and in the hidden constraints $(e_{hc})$ appear to be within the round-off error level, whatever stepsize is used.

\begin{figure}[p]
\centerline{\includegraphics[width=12cm,height=8cm]{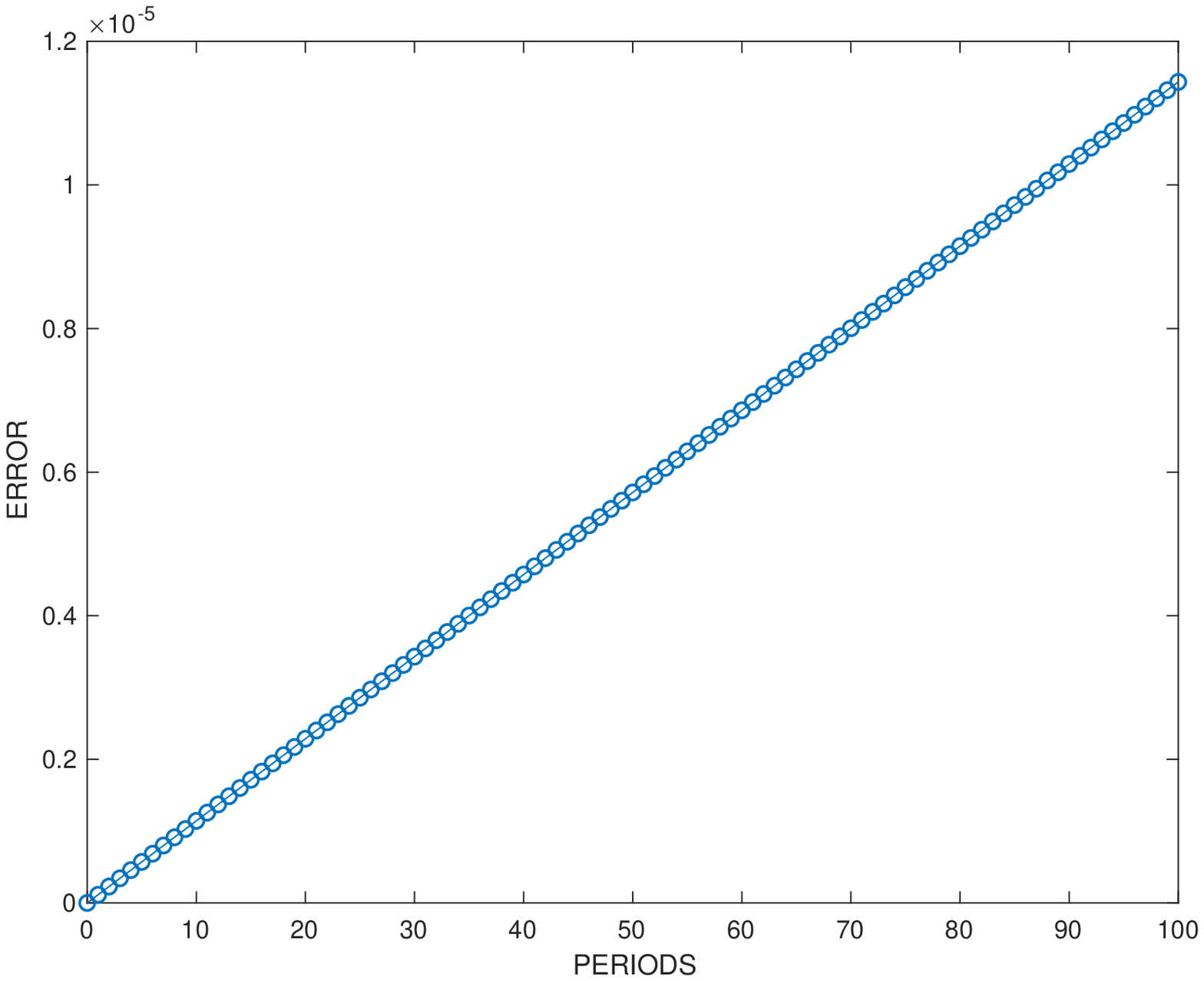}}
\caption{Conical pendulum (\ref{cpend_H})--(\ref{lambda0}). Linear growth of the solution error. The solution was computed over 100 periods using HBVM(2,2) and the stepsize $h=T/100\approx 0.053$. }
\label{fig1}

\bigskip
\centerline{\includegraphics[width=12cm,height=8cm]{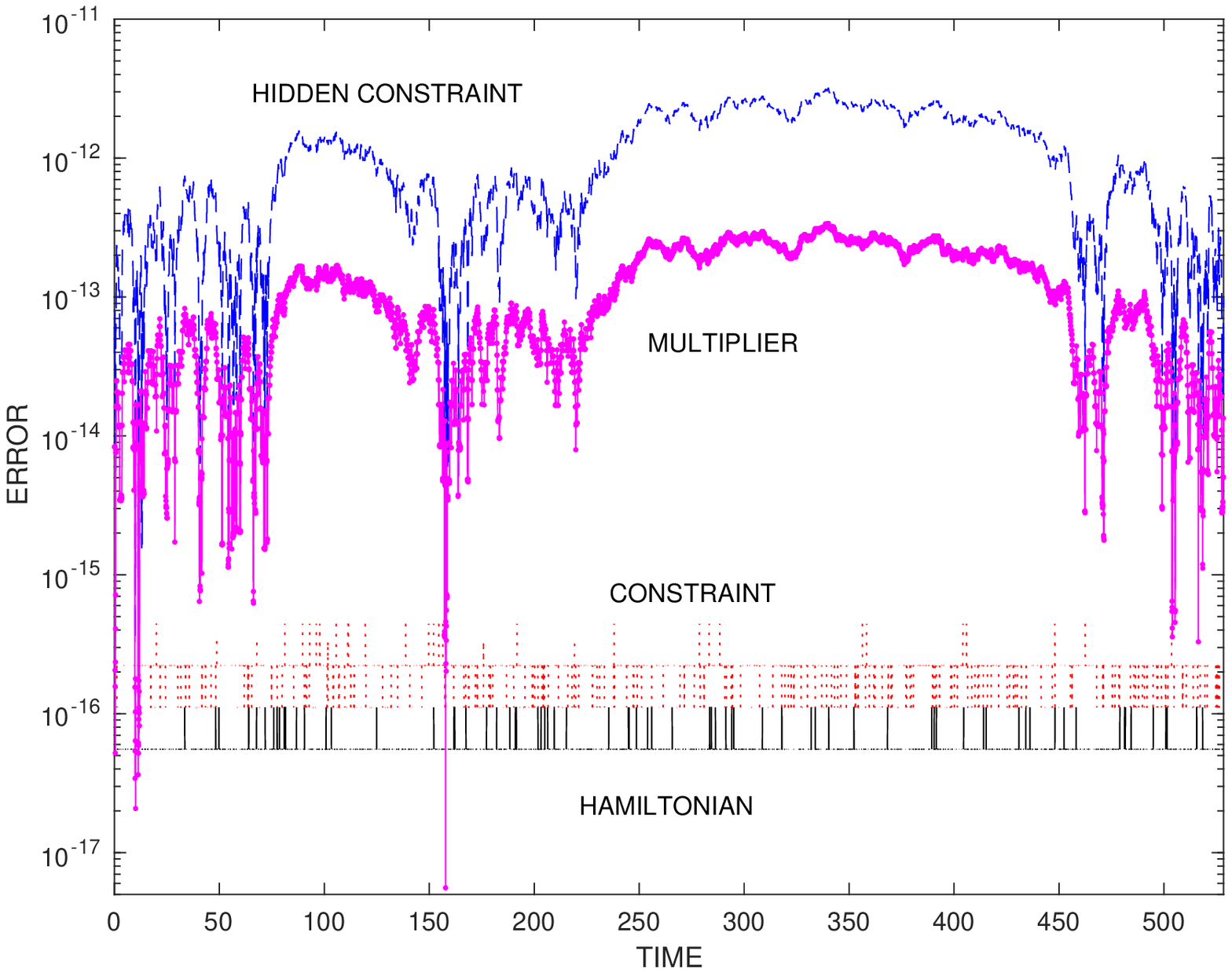}}
\caption{Conical pendulum (\ref{cpend_H})--(\ref{lambda0}). Multiplier error (dashed line), hidden constraint error (dash-dotted line), Hamiltonian error (solid line), and constraint error (dotted line) versus time. The simulation was  was carried our using HBVM(2,2) with the stepsize $h=T/100\approx 0.053$.}
\label{fig2}
\end{figure}

In Figure~\ref{fig1},  we plot the solution error from the computation over 100 periods using HBVM(2,2) with the stepsize $h=T/100\approx 0.053$, in Figure~\ref{fig2}, the errors of the multiplier, Hamiltonian, constraint, and hidden constraints. One can see the linear growth of the solution error. The errors of the multiplier, constraint, Hamiltonian, and hidden constraints are negligible.

\begin{table}[p]
\caption{\label{tab1} Conical pendulum (\ref{cpend_H})--(\ref{lambda0}). Errors from HBVM$(s,s)$ method for $s=1,2,3,4$, when solving the problem over 10 periods with stepsizes $h=T/n$.}  % see erroconic.m
\smallskip
\centerline{\begin{tabular}{|r|cccccc|}
\hline
\hline
       \multicolumn{7}{|c|}{$s=1$}\\
       \hline
$n$ & $e_s$       & rate  & $e_\lambda$ & $e_H$ & $e_g$ & $e_{hc}$ \\
\hline
 10 & 1.1543e\,00 & -- & 6.5503e-15 & 1.1102e-16 & 1.5543e-15 & 6.9435e-15\\
 20 & 4.0996e-01 & 1.49 & 1.2212e-14 & 1.1102e-16 & 6.6613e-16 & 7.3344e-15\\
 30 & 1.9021e-01 & 1.89 & 7.2831e-14 & 1.1102e-16 & 8.8818e-16 & 2.6870e-14\\
 40 & 1.0794e-01 & 1.97 & 2.3959e-13 & 1.1102e-16 & 6.6613e-16 & 6.8291e-14\\
 50 & 6.9285e-02 & 1.99 & 2.6112e-13 & 1.1102e-16 & 4.4409e-16 & 5.5241e-14\\
 60 & 4.8178e-02 & 1.99 & 1.8818e-13 & 1.1102e-16 & 1.2212e-15 & 4.0510e-14\\
 70 & 3.5420e-02 & 2.00 & 6.1351e-13 & 1.1102e-16 & 6.6613e-16 & 9.4355e-14\\
 80 & 2.7130e-02 & 2.00 & 8.1246e-13 & 1.1102e-16 & 5.5511e-16 & 1.0761e-13\\
 90 & 2.1441e-02 & 2.00 & 7.5329e-13 & 1.1102e-16 & 5.5511e-16 & 8.8662e-14\\
100 & 1.7371e-02 & 2.00 & 1.4311e-12 & 1.1102e-16 & 5.5511e-16 & 1.5112e-13\\
\hline
\hline
       \multicolumn{7}{|c|}{$s=2$}\\
       \hline
$n$ & $e_s$       & rate  & $e_\lambda$ & $e_H$ & $e_g$ & $e_{hc}$ \\
\hline
  10 & 1.1168e-02 & -- & 6.5503e-15 & 1.1102e-16 & 6.6613e-16 & 7.7539e-15\\
  20 & 7.1061e-04 & 3.97 & 1.8208e-14 & 5.5511e-17 & 3.3307e-16 & 1.1374e-14\\
  30 & 1.4083e-04 & 3.99 & 6.4060e-14 & 1.1102e-16 & 3.3307e-16 & 2.2125e-14\\
  40 & 4.4610e-05 & 4.00 & 2.2171e-13 & 1.1102e-16 & 4.4409e-16 & 6.2087e-14\\
  50 & 1.8282e-05 & 4.00 & 9.2704e-14 & 1.1102e-16 & 4.4409e-16 & 2.1613e-14\\
  60 & 8.8190e-06 & 4.00 & 3.6859e-13 & 1.1102e-16 & 4.4409e-16 & 6.4763e-14\\
  70 & 4.7611e-06 & 4.00 & 8.5076e-13 & 1.1102e-16 & 4.4409e-16 & 1.3128e-13\\
  80 & 2.7912e-06 & 4.00 & 1.2552e-12 & 1.1102e-16 & 2.2204e-16 & 1.6921e-13\\
\hline
\hline
       \multicolumn{7}{|c|}{$s=3$}\\
       \hline
$n$ & $e_s$       & rate  & $e_\lambda$ & $e_H$ & $e_g$ & $e_{hc}$ \\
\hline
 10 & 3.1758e-05 & -- & 6.1062e-15 & 1.1102e-16 & 1.5543e-15 & 9.8364e-15\\
 20 & 5.0199e-07 & 5.98 & 1.7319e-14 & 1.1102e-16 & 6.6613e-16 & 1.0819e-14\\
 30 & 4.4164e-08 & 5.99 & 9.8921e-14 & 1.1102e-16 & 2.2204e-16 & 3.4922e-14\\
 40 & 7.8663e-09 & 6.00 & 2.0650e-13 & 1.1102e-16 & 4.4409e-16 & 5.4473e-14\\
 50 & 2.0628e-09 & 6.00 & 2.3292e-13 & 1.1102e-16 & 4.4409e-16 & 4.9311e-14\\
 60 & 6.9103e-10 & 6.00 & 4.7151e-13 & 1.1102e-16 & 4.4409e-16 & 8.3119e-14\\
\hline
\hline
       \multicolumn{7}{|c|}{$s=4$}\\
       \hline
$n$ & $e_s$       & rate  & $e_\lambda$ & $e_H$ & $e_g$ & $e_{hc}$ \\
\hline
 10 & 4.9944e-08 & --   & 1.1768e-14 & 1.1102e-16 & 1.5543e-15 & 1.5603e-14\\
 20 & 1.9676e-10 & 7.99 & 1.7431e-14 & 1.1102e-16 & 6.6613e-16 & 1.1910e-14\\
 30 & 7.6630e-12 & 8.00 & 5.7399e-14 & 1.1102e-16 & 3.3307e-16 & 2.1564e-14\\
 40 & 7.3944e-13 & 8.13 & 4.1411e-14 & 1.1102e-16 & 4.4409e-16 & 1.0999e-14\\
\hline
\hline
\end{tabular}}
\end{table}

\subsection{Modified pendulum} % conicalpendN
We now consider a modified version of the previous problem. With this simulation, we aim at exploiting the conservation of energy and constraints using a suitable high-order quadrature rule (\ref{cgjd}). More precisely, we consider the following non-quadratic polynomial Hamiltonian:
\begin{equation}\label{cpend_H1}
H(q,p) := \frac{1}2p^\top p + (e_3^\top q)^4, \qquad q,p\in\RR^3,
\end{equation}
with the non-quadratic polynomial constraint
\begin{equation}\label{cpend_c1}
g(q) := \sum_{i=1}^3 (e_i^\top q)^{2(4-i)} -0.625 = 0.
\end{equation}
Here, we use the same initial condition as in (\ref{cpend_y0}) but the vector of the multipliers is no more constant, so that the order of the method reduces to $2$ (and $1$ for the vector of the multipliers).
Moreover, since the constraints and the Hamiltonian are polynomials of degree not larger than $6$, any HBVM$(3s,s)$ method will conserve both quantities. This is confirmed by Table~\ref{tab2}, where the results for HBVM(3,1), HBVM(6,2), and HBVM(9,3) are listed. The interval of integration was again $[0,10]$. Table~\ref{tab2} contains the errors in
\begin{itemize}
\item the solution ($e_s$),
\item the multiplier ($e_\lambda$),
\item the Hamiltonian $(e_H)$,
\item the constraints $(e_g)$;
\item the {\em hidden constraints} ($e_{hc}$), formally defined via (\ref{ehc}).
\end{itemize}
Again, as predicted in Theorem~\ref{capdcazz4}, we see that 
\begin{itemize}
\item all methods are second-order accurate in the state variables, with HBVM(3,1) less accu\-ra\-te than the others;
\item all methods are first-order accurate in the Lagrange multiplier;
\item all methods exactly conserve the Hamiltonian and the constraints;
\item all methods are second-order accurate in the hidden constraints.
\end{itemize}

\begin{table}[p]
\caption{\label{tab2} Modified pendulum (\ref{cpend_H1})--(\ref{cpend_c1}). Errors from HBVM$(3s,s)$ method for $s=1,2,3$, when solving the problem over the interval $[0,10]$ with stepsizes $h=10^{-1}2^{-n}$.}  % see erroconic.m
\smallskip
\centerline{\begin{tabular}{|r|rrrrrrrr|}
\hline
\hline
       \multicolumn{9}{|c|}{$s=1$}\\
       \hline
$n$ & $e_s$  & rate  & $e_\lambda$ & rate  & $e_H$ & $e_g$ & $e_{hc}$ & rate   \\
\hline
 0 & 2.0539e-02 &  --  & 1.0864e-01 &  --  & 1.1102e-16 & 1.6431e-14 & 1.5279e-02 &  -- \\
 1 & 4.9675e-03 & 2.05 & 6.4279e-02 & 0.76 & 2.2204e-16 & 7.5495e-15 & 3.9290e-03 & 1.96\\
 2 & 1.2365e-03 & 2.01 & 3.5621e-02 & 0.85 & 2.2204e-16 & 3.6637e-15 & 9.7072e-04 & 2.02\\
 3 & 3.0878e-04 & 2.00 & 1.8727e-02 & 0.93 & 2.2204e-16 & 2.1094e-15 & 2.4193e-04 & 2.00\\
 4 & 7.7173e-05 & 2.00 & 9.5965e-03 & 0.96 & 2.2204e-16 & 8.8818e-16 & 6.0436e-05 & 2.00\\
 5 & 1.9292e-05 & 2.00 & 4.8569e-03 & 0.98 & 2.2204e-16 & 5.5511e-16 & 1.5106e-05 & 2.00\\
 6 & 4.8229e-06 & 2.00 & 2.4432e-03 & 0.99 & 2.2204e-16 & 6.6613e-16 & 3.7764e-06 & 2.00\\
 7 & 1.2057e-06 & 2.00 & 1.2251e-03 & 1.00 & 2.2204e-16 & 6.6613e-16 & 9.4417e-07 & 2.00\\
 8 & 3.0139e-07 & 2.00 & 6.1472e-04 & 1.00 & 2.2204e-16 & 6.6613e-16 & 2.3608e-07 & 2.00\\
\hline
\hline
       \multicolumn{9}{|c|}{$s=2$}\\
       \hline
$n$ & $e_s$  & rate  & $e_\lambda$ & rate  & $e_H$ & $e_g$ & $e_{hc}$ & rate   \\
\hline
 0 & 6.0495e-03 &  --  & 1.5224e-01 &  --  & 1.1102e-16 & 5.7732e-15 & 1.7516e-02 &  -- \\
 1 & 1.4027e-03 & 2.11 & 7.6627e-02 & 0.99 & 2.2204e-16 & 3.8858e-15 & 4.6710e-03 & 1.91\\
 2 & 3.4600e-04 & 2.02 & 3.8806e-02 & 0.98 & 1.1102e-16 & 2.3315e-15 & 1.1666e-03 & 2.00\\
 3 & 8.6197e-05 & 2.01 & 1.9532e-02 & 0.99 & 2.2204e-16 & 1.6653e-15 & 2.9091e-04 & 2.00\\
 4 & 2.1530e-05 & 2.00 & 9.7988e-03 & 1.00 & 2.2204e-16 & 5.5511e-16 & 7.2716e-05 & 2.00\\
 5 & 5.3813e-06 & 2.00 & 4.9076e-03 & 1.00 & 2.2204e-16 & 5.5511e-16 & 1.8175e-05 & 2.00\\
 6 & 1.3453e-06 & 2.00 & 2.4559e-03 & 1.00 & 2.2204e-16 & 6.6613e-16 & 4.5440e-06 & 2.00\\
 7 & 3.3632e-07 & 2.00 & 1.2285e-03 & 1.00 & 2.2204e-16 & 6.6613e-16 & 1.1360e-06 & 2.00\\
 8 & 8.4002e-08 & 2.00 & 6.1386e-04 & 1.00 & 2.2204e-16 & 6.6613e-16 & 2.8414e-07 & 2.00\\
\hline
\hline
       \multicolumn{9}{|c|}{$s=3$}\\
       \hline
$n$ & $e_s$  & rate  & $e_\lambda$ & rate  & $e_H$ & $e_g$ & $e_{hc}$ & rate   \\
\hline
 0 & 6.0698e-03 &  --  & 1.5231e-01 &  --  & 2.2204e-16 & 4.6629e-15 & 1.7532e-02 &  -- \\
 1 & 1.4040e-03 & 2.11 & 7.6632e-02 & 0.99 & 1.1102e-16 & 3.9968e-15 & 4.6715e-03 & 1.91\\
 2 & 3.4608e-04 & 2.02 & 3.8806e-02 & 0.98 & 1.1102e-16 & 1.7764e-15 & 1.1666e-03 & 2.00\\
 3 & 8.6202e-05 & 2.01 & 1.9532e-02 & 0.99 & 2.2204e-16 & 1.3323e-15 & 2.9091e-04 & 2.00\\
 4 & 2.1530e-05 & 2.00 & 9.7988e-03 & 1.00 & 2.2204e-16 & 5.5511e-16 & 7.2716e-05 & 2.00\\
 5 & 5.3814e-06 & 2.00 & 4.9076e-03 & 1.00 & 2.2204e-16 & 6.6613e-16 & 1.8175e-05 & 2.00\\
 6 & 1.3453e-06 & 2.00 & 2.4560e-03 & 1.00 & 2.2204e-16 & 6.6613e-16 & 4.5439e-06 & 2.00\\
 7 & 3.3623e-07 & 2.00 & 1.2283e-03 & 1.00 & 2.2204e-16 & 6.6613e-16 & 1.1361e-06 & 2.00\\
 8 & 8.4116e-08 & 2.00 & 6.1452e-04 & 1.00 & 2.2204e-16 & 6.6613e-16 & 2.8410e-07 & 2.00\\
\hline
\hline
\end{tabular}}
\end{table}

\subsection{Tethered satellites system} % tether.m
\begin{figure}[p]
\centerline{\includegraphics[width=12cm,height=12cm]{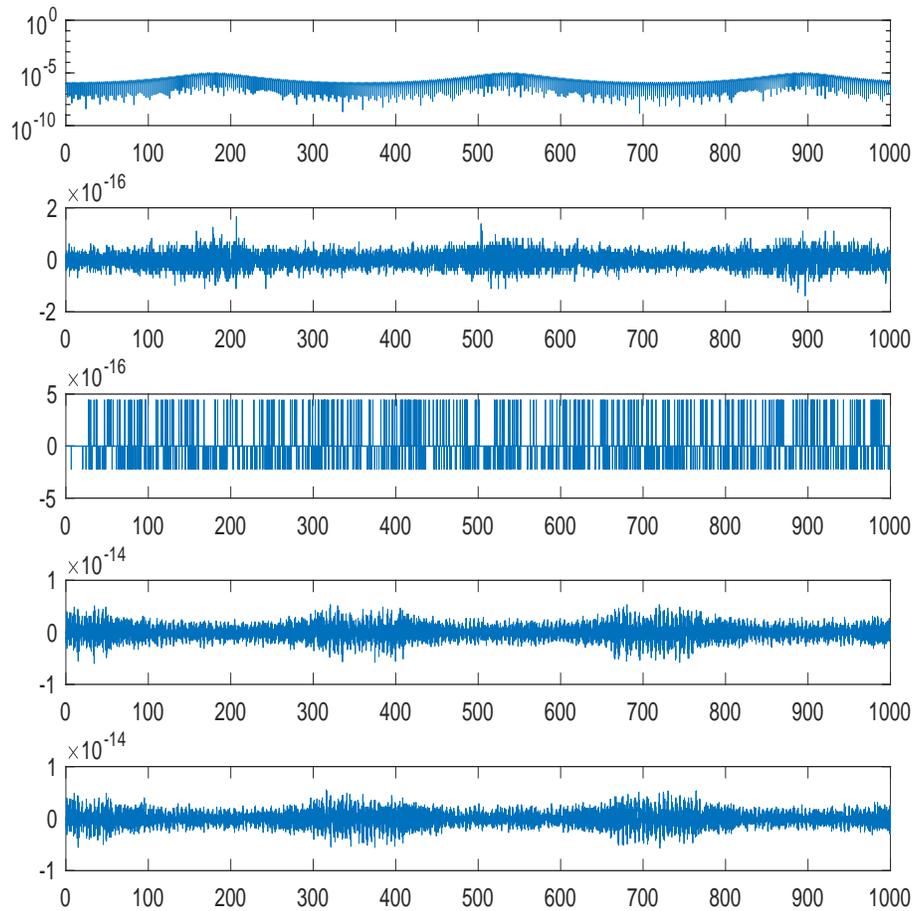}}
\caption{Tethered satellites system (\ref{gts})--(\ref{p1230}). The problem was solved using HBVM(6,2) with the stepsize $h=0.1$ over $10^4$ steps. From top to bottom: hidden constraints error (first plot), Hamiltonian error (second plot), constraints errors (third to fifth plots).}
\label{fig3}
\end{figure}

Finally, we discuss a closed-loop rotating triangular tethered satellites system,\footnote{This example can be found in \cite{CLZ2015,WDLW2016}.} including three satellites (considered as mass-points) of masses $m_i$, $i=1,2,3$, joined by inextensible, tight, and massless tethers, of lengths $L_i$, $i=1,2,3$, respectively, which orbit around a massive body.\footnote{The Earth, in the original example.} As before, for sake of simplicity, we assume unit masses and lengths, and normalize the gravity constant. Consequently, if $q_i:=(x_i,y_i,z_i)^\top\in\RR^3$, $i=1,2,3$, are the positions of the three satellites, the constraints are given by
\begin{equation}\label{gts}
g(q) := \pmatrix{c} (q_1-q_2)^\top(q_1-q_2) - 1 \\
(q_2-q_3)^\top(q_2-q_3) - 1 \\
(q_3-q_1)^\top(q_3-q_1) - 1 \endpmatrix= 0\in\RR^3,
\end{equation}
and the Hamiltonian is specified by 
\begin{equation}\label{Hts}
H(q,p) = \sum_{i=1}^3 \left(\frac{1}2  p_i^\top p_i - \frac{1}{\sqrt{q_i^\top q_i}}\right).
\end{equation}
The consistent initial conditions have the form
\begin{equation}\label{q1230}
q_1(0) = \pmatrix{r}0 \\ \frac{1}2\\       z_0  \endpmatrix, \qquad
q_2(0) = \pmatrix{r}0 \\ -\frac{1}2\\       z_0  \endpmatrix, \qquad
q_3(0) = \pmatrix{c}0 \\ 0 \\       z_0-\frac{\sqrt{3}}2  \endpmatrix,
\end{equation}
and
\begin{equation}\label{p1230}
p_1(0)=p_2(0)= \pmatrix{c} 0 \\ 0 \\    0  \endpmatrix,\qquad
p_3(0) = \pmatrix{c} v_0 \\ 0 \\    0  \endpmatrix,
\end{equation}
where $z_0=20$ and $v_0$ is such that the initial Hamiltonian is zero. This provides a configuration in which the first two satellites remain parallel to each other, moving in the planes $y=\frac{1}2$ and $y=-\frac{1}2$, respectively, and the third one moves around the tether joining the first two, in the plane $y=0$. In such a case, the Hamiltonian is non-polynomial. Nevertheless, using the HBVM(6,2) method with the stepsize $h=0.1$ over $10^4$ steps, we obtain a qualitatively correct solution which conserves the Hamiltonian and the constraints within the round-off error level, see Figure~\ref{fig3}. Here, we also plot the hidden constraints errors $ \| \nabla g(q_n)^\top M^{-1}p_n\|$.
At last, in Table~\ref{tab3}, we list the following errors arising when solving the problem with HBVM$(6,s)$ methods for $s=1,2,3$ and the stepsizes $h=10^{-1}2^{-n}$, over the interval $[0,10]$:
\begin{itemize}
\item the solution error ($e_s$),
\item the multipliers error ($e_\lambda$),
\item the Hamiltonian error $(e_H)$,
\item the constraints error $(e_g)$;
\item the {\em hidden constraints} errors ($e_{hc}$), formally defined by (\ref{ehc}).
\end{itemize}
Again, as shown in Theorem~\ref{capdcazz4}, 
\begin{itemize}
\item all methods are second-order accurate in the state variables, with HBVM(6,1) much less accu\-ra\-te than the other two (which are almost equivalent);
\item all methods are first-order accurate in the Lagrange multipliers;
\item all methods exactly conserve the Hamiltonian and the constraints;
\item all methods are second-order accurate in the hidden constraints.
\end{itemize}
%\begin{rem}

To draw a general conclusion: it seems that using of HBVM$(k,s)$, with $s>1$, in context of the numerical solution of the Hamiltonian problems with holonomic constraints can be recommended, although the method is only second-order accurate.
%\end{rem}

\begin{table}[p]
\caption{\label{tab3} Tethered satellite system (\ref{gts})--(\ref{p1230}). Errors from the HBVM$(6,s)$ method for $s=1,2,3$, when solving the problem over the interval [0,10] with stepsizes $h=10^{-1}2^{-n}$.}  % see erroconic.m
\smallskip
\centerline{\begin{tabular}{|r|rrrrrrrr|}
\hline
\hline
       \multicolumn{9}{|c|}{$s=1$}\\
       \hline
$n$ & $e_s$  & rate  & $e_\lambda$ & rate  & $e_H$ & $e_g$ & $e_{hc}$ & rate   \\
\hline
0 & 9.2893e-04 &  --  & 2.1218e-06 &  --  & 5.5511e-17 & 1.2212e-14 & 9.6503e-07 &  -- \\
1 & 2.3234e-04 & 2.00 & 1.0689e-06 & 0.99 & 5.5511e-17 & 1.2212e-14 & 2.4108e-07 & 2.00\\
2 & 5.8093e-05 & 2.00 & 5.3623e-07 & 1.00 & 6.9389e-17 & 1.5765e-14 & 6.0272e-08 & 2.00\\
3 & 1.4524e-05 & 2.00 & 2.6834e-07 & 1.00 & 6.9389e-17 & 1.5099e-14 & 1.5090e-08 & 2.00\\
\hline
\hline
       \multicolumn{9}{|c|}{$s=2$}\\
       \hline
$n$ & $e_s$  & rate  & $e_\lambda$ & rate  & $e_H$ & $e_g$ & $e_{hc}$ & rate   \\
\hline
0 & 1.8586e-07 &  --  & 2.1635e-06 &  --  & 6.2450e-17 & 1.1990e-14 & 1.3053e-06 &  -- \\
1 & 3.9859e-08 & 2.22 & 1.0812e-06 & 1.00 & 4.8572e-17 & 1.4877e-14 & 3.2584e-07 & 2.00\\
2 & 9.5501e-09 & 2.06 & 5.4039e-07 & 1.00 & 6.9389e-17 & 1.4433e-14 & 8.1466e-08 & 2.00\\
3 & 2.3636e-09 & 2.01 & 2.7052e-07 & 1.00 & 6.9389e-17 & 1.4655e-14 & 2.0374e-08 & 2.00\\
\hline
\hline
       \multicolumn{9}{|c|}{$s=3$}\\
       \hline
$n$ & $e_s$  & rate  & $e_\lambda$ & rate  & $e_H$ & $e_g$ & $e_{hc}$ & rate   \\
\hline
0 & 1.5089e-07 & --      &  2.1635e-06  & --      &   5.5511e-17 &  1.3767e-14 & 1.3053e-06  & --\\
1 & 3.7674e-08 & 2.00 &  1.0813e-06  & 1.00 &   6.2450e-17 &  1.4211e-14 &  3.2585e-07 &  2.00\\
2 & 9.4170e-09 & 2.00 &   5.4061e-07 & 1.00 &   6.9389e-17 &  1.3323e-14 &  8.1471e-08 & 2.00\\
3 & 2.3630e-09 & 2.00 &   2.7140e-07 & 1.00 &   6.9389e-17 &  1.5099e-14 &  2.0381e-08 & 2.00\\
\hline
\hline
\end{tabular}}
\end{table}

\section{Conclusions}\label{end}
In this paper, we have considered the numerical solution of Hamiltonian problems with holonomic constraints, by resorting to a line-integral formulation of the conservation of the constraints. This approach enables to derive an expression for the Lagrange multipliers in which {\em second derivatives are not used}. From the discretization of the resulting formulae, we have obtained a suitable variant of the Hamiltonian Boundary Value Methods (HBVMs), formerly designed as an energy-conserving Runge-Kutta methods for unconstrained Hamiltonian problems. Numerical experiments supporting the theoretical findings are enclosed. \\

This paper has been initiated in spring 2017, during the visit of the first author at the Institute for Analysis and Scientific Computing, Vienna University of Technology, Vienna, Austria.

\end{document}